\documentstyle{amsppt}
\magnification=\magstep1

\topmatter
\title Differential Calculi over Quantum Groups and Twisted Cyclic Cocycles \endtitle
\author J. Kustermans, G.J. Murphy and L. Tuset \endauthor
\address Department of Mathematics, KU Leuven, Belgium \endaddress
\address
Department of Mathematics, National University of Ireland, Cork, Ireland
\endaddress
\address  Faculty of Engineering, University College, Oslo, Norway. \endaddress
\abstract We study some aspects of the theory of non-commutative differential calculi over complex algebras, especially
over the Hopf algebras associated to compact quantum groups in the sense of S.L.~Woronowicz. Our principal emphasis is
on the theory of twisted graded traces and their associated twisted cyclic cocycles. One of our principal results is a
new method of constructing differential calculi, using twisted graded traces. \endabstract \keywords Hopf algebra;
differential calculus; twisted graded trace \endkeywords \subjclass 46L, 81R50 \endsubjclass
\endtopmatter

\NoRunningHeads \overfullrule=0pt


\def\a{\alpha}  \def\g{\gamma} \def\d{\delta}
\def\e{\varepsilon} \def\f{\varphi} \def\k{\kappa} \def\l{\lambda}
 \def\r{\rho} \def\s{\sigma} \def\ta{\theta} \def\t{\tau}
\def\om{\omega} \def\Om{\Omega}  \def\D{\Delta}

\def\bC{{\bold C}} \def\bN{{\bold N}} 
   \def\bZ{{\bold Z}}

 \def\id{\operatorname{id}}

\def\bZ{\bold Z} \def\bN{\bold N} 
 \def\bC{\bold C}

\def\sa{\Cal A}  \def\ominv{\Om^{\text{inv}}}

\def\omd{{(\Om,d)}} \def\bom{{\bar \Om}} \def\bomd{{(\bar\Om,d)}}
\def\tom{{\tilde{\Om}}}\def\tomd{{(\tilde{\Om},d)}}
\def\sad{{(\sa,\D)}} \def\dom{{\D_\Om}}
\def\idts{{\id_\sa\otimes \int}} \def\htid{{h\otimes \id_\Om}} \def\bint{\int'}
\def\idth{{\id_\sa\otimes h}}


\def\mapp#1{\buildrel{\tsize #1}\over\longrightarrow}
\def\triangle#1#2#3#4#5#6{\matrix #1&\mapp #2&#3\\
\vspace{2\jot}                                \quad {}_{#6}\nwarrow &&\swarrow {}_{#4} \quad \\
                                  &#5\endmatrix}

\def\de{\Delta}
\def\ot{\otimes}

\def\al{\alpha}

\def\om{\omega}
\def\Om{\Omega}

\def\ABE{1}
\def\CON{2}
\def\KAS{3}
\def\KS{4}
\def\KT{5}
\def\MT{6}
\def\WORA{7}
\def\WORB{8}
\def\WORC{9}

\document

\subheading{Introduction} A compact group is a compact space with a continuous multiplication satisfying certain extra conditions. In the theory of compact quantum groups developed by S.L.~Woronowicz~\cite{\KAS, \KS, \KT, \MT, \WORB}, one replaces the compact space by a unital C*-algebra~$A$ that is in general non-commutative, and replaces the group multiplication by a co-multiplication on~$A$ satisfying certain cancelation conditions. Contained in~$A$ is a dense $*$-subalgebra~$\sa$, the {\it representation} algebra, that is a Hopf algebra under the restriction co-multiplication. Both $A$ and $\sa$ admit a Haar integral and this is vital for many aspects of the theory we develop in this paper.

The considerations in this paper are motivated by the theory of compact quantum groups, but it is not these objects that we study here; rather, we study differential calculi over such groups. Our context is therefore non-commutative differential geometry in the spirit of that subject as developed by Alain Connes~\cite{\CON}. The study of differential calculi in the quantum group setting was initiated by Woronowicz---indeed, he constructed the first example of such a calculus~\cite{\WORA}. However, it was immediately apparent in his work that Connes' theory of non-commutative geometry does not cover the calculi occurring in the quantum setting. To explain briefly what is involved, recall that although the algebra of forms in the classical setting of differential manifolds is not commutative, it is ``nearly'' so, in the sense that $\om\om'= (-1)^{kl}\om'\om$, if $\om$ and $\om'$ are a $k$-form and an $l$-form, respectively. In Connes' non-commutative geometry, it is no longer true that $\om\om'= (-1)^{kl}\om'\om$. However, for a graded trace (this is an appropriate kind of ``integral'' on the ``non-commutative manifold''), we have
${\int \om\om'}={(-1)^{kl}\int \om'\om}$, where $\om$ and $\om'$ are a $k$-form and an $l$-form, respectively. This integral condition is of fundamental importance in the cyclic cocyle theory developed so successfully by Connes in the past two decades.
However, even this weaker commutativity condition does not hold in the context of differential geometry over quantum groups. If one thinks of a graded trace as the analogue of a trace on a C*-algebra, then one can explain the situation in the quantum setting by saying that one must replace a trace by a KMS~state. More precisely, in this setting there is an automorphism $\s$ of degree zero of the algebra of forms such that
${\int \om\om'}={(-1)^{kl}\int \s(\om')\om}$, where $\om$ and $\om'$ are a $k$-form and an $l$-form, respectively. This is, of course, analogous to the situation with a KMS~state~$h$ on a C*-algebra, where one has an automorphism $\s$ on a dense $*$-subalgebra for which $h(ab)=h(\s(b)a)$, for all elements $a$ and $b$ in the subalgebra.

In his seminal paper on differential calculi over quantum groups~\cite{\WORA}, Woronowicz remarks that the integral he defines on his 3-dimensional calculus over the quantum group $SU_q(2)$ does not fit into the framework of Connes' non-commutative geometry, but he does not develop this observation. In this paper we introduce the concept of a twisted graded trace (the analogue of a KMS state) to replace Connes' graded traces. It is then necessary to develop a theory of twisted cyclic cocycles and we do this here. One of our principal results is a new method of constructing differential calculi; in essence, in this approach we start with a twisted graded trace and construct a calculus (in Woronowicz's approach one goes in the opposite direction). We feel that our approach may be more natural, since, to some extent, it involves giving a ``presentation'' of the calculus in terms of generators and relations.

We give a brief overview of the paper now. In Section~1 we introduce the basic terminology and prove two theorems that are very useful for constructing twisted graded traces. We also introduce a quotient construction for obtaining a differential calculus from a twisted graded trace. In Section~2 we introduce twisted cyclic cocycles and develop their relationship with twisted graded traces. In both this section and the next, we develop a theory of twisted cyclic cohomology. This contains Connes' theory as a special case, but, as we have indicated above,  the more general theory is necessary to deal with the examples that occur in the quantum group setting. However, the theory developed in Sections~1--3
is not restricted to the quantum group setting and applies in the more general context of differential calculi over arbitrary unital algebras. In Section~4 we develop aspects of the theory of left-invariant twisted graded traces over left-covariant differential calculi. In this situation the underlying algebra is assumed to be a Hopf algebra. An important result here is that the differential calculus constructed from a left-invariant twisted graded trace on the universal calculus is shown to be itself left-covariant. Also, we give a characterization of the twisted cyclic cocycles that correspond to left-invariant twisted graded traces. In the final section, Section~5, we show in detail how our ideas can be used to give an alternative construction of Woronowicz's first, 3-dimensional, differential calculus over quantum $SU(2)$. This is a construction that may be of prime importance in the future for obtaining new examples of differential calculi.

\subheading{1. Differential calculi}

In this section we set up the basic terminology for studying differential calculi over algebras that are not necessarily commutative. One can think of this as the study of differential forms in the setting of quantum spaces or manifolds. We give a general procedure for constructing such calculi.
We begin by recalling some basic definitions.

Let $\Om$ be a (positively) graded algebra, $\Om={\oplus_{n=0}^\infty \Om_n}$. A {\it graded derivation} on $\Om$ is a linear map ${d\colon \Om\to \Om}$ for which
$d(\om'\om)=d(\om')\om+(-1)^n\om' d\om$,
for all $\om'\in \Om_n$ and all $\om\in \Om$.

A {\it graded differential algebra} is a pair $(\Om,d)$, where $\Om$ is a graded algebra, $d$ is a graded derivation on $\Om$ of degree 1 (as a linear map) and $d^2=0$. The elements of $\Om$ are referred to as the {\it forms} of $(\Om,d)$ and the elements of $\Om_n$ as the
{\it $n$-forms}. The operator $d$ is referred to as the {\it differential}.

Now suppose that $\sa$ is an arbitrary associative unital algebra.
Then there is a graded differential algebra $\bomd$, for which
$\bar\Om_0=\sa$, that has the following universal property: If
$\s$ is an  algebra homomorphism from $\sa$ into the algebra
$\Om_0$ of $0$-forms of a graded differential algebra $\omd$, then
there exists a unique algebra homomorphism $\bar\s$ from $\bom$ to
$\Om$ extending $\s$ such that $\bar\s d=d\bar\s$. This property
uniquely determines $\bomd$ (up to isomorphism). Note that
$\bar\s$ is clearly necessarily of grade zero. We shall usually
denote the extension $\bar\s$ by the same symbol $\s$ as the
original homomorphism.

We shall use the following two useful properties of $\bomd$:

(1) Let $n \geq 1$. Then every element of $\bom_n$ is a sum of elements of the form $a_0da_1\cdots da_n$,  and $da_1\cdots da_n$, where the elements $a_0,a_1,\dots,a_n$ belong to $\sa$;

(2) Let $n$ be a positive integer and $T_1$ a multilinear map from ${\sa^{n+1}}$ to a linear space $Y$ and $T_2$ a linear map from ${\sa^n}$ to the same linear space $Y$. Then there is a unique linear map $\hat T$ from $\bom_n$ to $Y$ for which $\hat T(a_0da_1\cdots da_n)=T_1(a_0,a_1,\dots,a_n)$ and $\hat T(da_1\cdots da_n)=T_2(a_1,\dots,a_n)$, for all $a_0,a_1,\dots,a_n\in \sa$.

In practice, the universal graded differential algebra $\bomd$ is too big to be useful. However, it can be used to construct smaller, finite-dimensional differential algebras that are useful.

A {\it differential calculus} over $\sa$ is a graded differential algebra $\omd$ for which

(1) $\Om_0=\sa$;

(2) Let $n \geq 1$. Then every element of $\Om_n$ is a sum of elements of the form $a_0da_1\cdots da_n$ and $da_1\cdots da_n$, where the elements $a_0,a_1,\dots,a_n$ belong to $\sa$.

If the differential calculus $\Om$ is unital (as an algebra), then the unit of $\Om$ has to belong to $\Om_0 = \sa$ and therefore has to be equal to the unit $1$ of $\sa$.

We shall say the differential calculus $\omd$ is {\it finite-dimensional}, of {\it dimension} $N$, if $\Om_N\ne 0$ and $\Om_n=0$ for $n>N$.

The universal graded differential algebra is clearly a differential calculus over $\sa$, but it is, equally clearly, not finite-dimensional, nor unital.

We now describe a general procedure for obtaining a new, ``smaller'' calculus from a given calculus.
Let $N$ be a positive integer and let $\omd$ be a differential calculus over~$\sa$ that is either not finite-dimensional, or is of finite dimension greater than~$N$. We define a new differential calculus $(\Om',d')$ of dimension~$N$ by setting $\Om'_k=\Om_k$, if $k\le N$ and $\Om'_k=0$, if $k>N$. We define the multiplication $\cdot$ in $\Om'$ by setting, for $\om_1\in \Om_k$ and $\om_2\in \Om_l$,
$\om_1\cdot \om_2=\om_1\om_2$, if $k+l\le N$, and by setting $\om_1\cdot \om_2=0$ if $k+l>N$. We set $d'(\om_1)=d(\om_1)$, if $k\le N$ and set $d'(\om_1)=0$, if $k>N$. We call $(\Om',d')$ the differential calculus of dimension~$N$ {\it obtained from $\omd$ by truncation}.

If $\omd$ is a differential calculus over $\sa$, we say a linear functional $\int$ on $\Om$ is {\it closed} if $\int d=0$.
If $\om_1,\dots,\om_M\in \Om$, then a simple induction shows that $d\om_1d\om_2\cdots d\om_M=d(\om_1d\om_2\cdots d\om_M)$. Hence, if $\int$ is closed, $\int d\om_1d\om_2\cdots d\om_M=0$. We shall frequently tacitly make use of this observation. If $\om$ is a $k$-form and $\om'$ an arbitrary form, then $\int (d\om)\om'=(-1)^{k+1}\int\om d\om'$, another result we shall use tacitly in the sequel. It follows from the fact that $d(\om\om')=(d\om)\om'+(-1)^k\om d\om'$ and $\int d(\om\om')=0$.

A linear functional $\int$ on $\Om$ is a {\it twisted graded trace} if
there is an algebra automorphism ${\s\colon \Om\to \Om}$ of degree zero for which $\s d=d\s$ and
$\int\om'\om=(-1)^{kl}\int\s(\om)\om'$,
for all non-negative integers $k$ and $l$ and for all $\om\in \Om_k$ and $\om'\in \Om_l$.

We say $\s$ is a {\it twist automorphism} associated to $\int$. It is useful to observe that $\int \s(\om)= \int \om$, for all $\om \in \Om$. To see this, observe first that $a = a 1$ and $da = d(a1) = (da)1 + a(d1)$ for all $a \in A$. It follows that any element of $\Om$ is a sum of products of two elements of $\Om$. Let $\om,\om' \in  \Om$. We may write ${\om=\sum_k \om_k}$ and ${\om'=\sum_k \om'_k}$, where $\om_k,\om'_k\in \Om_k$. Then ${\int \om\om'}={\sum_{k,l}\int \om_k\om'_l}=
{\sum_{k,l} (-1)^{kl}\int \s(\om'_l)\om_k}={\sum_{k,l}\int \s(\om_k)\s(\om'_l)}=
{\int \s(\om)\s(\om')}={\int \s(\om\om')}$.

\proclaim{Theorem 1.1} Let $\bomd$ be the universal calculus over a unital  algebra $\sa$. Suppose that $\int$ is a closed linear functional on $\bom$ and that ${\s_0\colon \sa\to \sa}$ is an algebra automorphism for which $\int \s_0(a)\om=\int \om a$, for all $a\in \sa$ and $\om\in \bom$. Then $\int$ is a twisted graded trace having a twist automorphism $\s$ that extends~$\s_0$. \endproclaim

\demo{Proof} The automorphism, ${\s_0\colon \bom_0\to \bom_0}$, extends uniquely to an automorphism,
${\s\colon \bom\to \bom}$, for which $\s d=d\s$, by the universal property of~$\bomd$.
We shall show that $\int$ is a twisted graded trace, with $\s$ as its twist automorphism.
Thus, to prove the theorem, we have only to show that, for each positive integer~$N$,
$$\int\om'\om=(-1)^{k(N-k)}\int\s(\om)\om', \tag 1$$
for all integers $k$ such that $0\le k \le N$, and for all $\om\in \bom_k$ and $\om'\in\bom_{N-k}$.
We shall prove this by induction on~$k$. It clearly holds for $k=0$ by hypothesis. Let's assume it holds for $k$ and we shall prove it for $k+1$, where we also suppose that $k+1\le N$. We first show that
$$\int \a d\om=(-1)^{(k+1)(N-k-1)}\int \s(d\om)\a, \tag 2$$
where $\om \in \bom_k$ and $\a\in \bom_{N-k-1}$. We suppose first that $k+1<N$. If $\al = d\om'$, where $\om'\in \bom_{N-k-2}$, the closedness of $\int$ implies that both sides of the above equation are 0 and hence equal. Since $\bom_{N-k-1}$ is the linear span of elements of the form $d\om'$ and $(d\om')a$, where $\om'\in \bom_{N-k-2}$ and $a\in \sa$, we may now clearly suppose that $\a=(d\om')a$. We have
$\int (d\om')ad\om=\int d\om'd(a\om)-\int(d\om')(da)\om=
-\int (d\om')(da)\om=
(-1)^{1+k(N-k)}\int \s(\om)(d\om')da$, by the inductive hypothesis. Since $d(\s(\om)\om')=(d\s(\om))\om'+(-1)^k\s(\om)d\om' = \s(d\om)\om'+(-1)^k\s(\om)d\om'    $, we get
$$\align
\int (&d\om')ad\om=(-1)^{1+k(N-k)}(-1)^k[\int d(\s(\om)\om')da-\int \s(d\om)\om'da]\\
&=(-1)^{1+k(N-k)}(-1)^{k+1}\int \s(d\om)\om'da\\
&=(-1)^{1+k(N-k)}(-1)^{k+1}(-1)^{N-k-2}[\int \s(d\om)d(\om'a)-\int \s(d\om)(d\om')a]\\&=
(-1)^{1+k(N-k)}(-1)^{k+1}(-1)^{N-k-1}\int \s(d\om)(d\om')a\\&=
(-1)^{(k+1)(N-k-1)}\int \s(d\om)(d\om')a.
\endalign$$
This shows that Equation~(2) holds, as required, when $k+1<N$. For
$k+1=N$ the argument is similar, but much simpler, and is
therefore omitted. It follows now from Equation~(2) that, for all
$a\in \sa$, we have
$$\align
\int &\a ad\om=(-1)^{(k+1)(N-k-1)}\int
\s(d\om)\a a\\
&= (-1)^{(k+1)(N-k-1)}\int \s_0(a)\s(d\om)\a=
\int \s(ad\om)\a.
\endalign$$ This
shows that Equation~(1) is satisfied for $k$ in place of $k+1$.
This completes our induction, so Equation~(1) is now seen to be
true for ${k=0,\dots,N}$. \qed\enddemo

We say that a linear functional $\int$ on $\Om$ is {\it left faithful} if, whenever $\om\in \Om$ is such that $\int \om'\om=0$, for all $\om'\in \Om$, we necessarily have $\om=0$.

\proclaim{Theorem 1.2} Suppose $\omd$ is a differential calculus over a unital algebra~$\sa$. Suppose that $\int$ is a left faithful, closed linear functional on $\Om$ and that ${\s_0\colon \sa\to \sa}$ is an algebra automorphism for which $\int \s_0(a)\om=\int \om a$, for all $a\in \sa$ and $\om\in \Om$. Then $\int$ is a twisted graded trace having a twist automorphism $\s$ that extends $\s_0$. \endproclaim

\demo{Proof} The automorphism, ${\s_0\colon \bom_0\to \bom_0}$,
extends uniquely to an automorphism, ${\bar\s\colon \bom\to
\bom}$, for which $\bar\s d=d\bar\s$, by the universal property of
the universal differential calculus $\bomd$. Likewise the
isomorphism, ${\id_\sa\colon \bom_0\to \Om_0}$, extends uniquely
to a surjective homomorphism, ${\pi\colon \bom\to \Om}$, such that
$\pi d =d\pi$. We define $\bint$ on $\bom$ by setting ${\bint
\om}={\int \pi(\om)}$, for all $\om \in \bom$. Clearly, $\bint$ is
a closed, linear functional on $\bom$ satisfying the hypothesis of
the preceding theorem. Hence, $\bint$ is a twisted graded trace,
with $\bar\s$ as its twist automorphism.

Suppose now that $\om\in\bom$ and $\pi(\om)=0$. We shall show that $\pi(\bar\s(\om))=0$. If $\om'\in \bom$, then ${\int \pi(\bar\s(\om'))\pi(\bar\s(\om))}=
{\bint \bar\s(\om'\om)}={\bint \om'\om}={\int\pi(\om')\pi(\om)=0}$, since $\pi(\om)=0$. It follows from faithfulness of $\int$ that $\pi(\bar\s(\om))=0$, as required.

We can now use this invariance of $\ker(\pi)$ under $\bar\s$ to
induce a homomorphism $\s$ on $\Om$ defined by setting
$\s(\pi(\om))=\pi(\bar\s(\om))$, for all $\om\in \bom$. It is
clear that ${\int\om'\om}={(-1)^{kl}\int\s(\om)\om'}$, for all
integers $k$ and~$l$ and for all $\om\in \Om_k$ and
$\om'\in\Om_l$. Clearly, since $\bar\s$ extends $\s_0$, so
does~$\s$. It is easily checked that $\s d = d\s$. Moreover, $\s$ is surjective, since
$\bar\s$ and $\pi$ are. Thus, to show that $\int$ is a twisted
graded trace with $\s$ as twist automorphism, we need only show
now that $\s$ is injective. To see this, suppose that $\om\in
\Om_k$ and $\s(\om)=0$. Then $\int
\om'\om=(-1)^{kl}\int\s(\om)\om'=0$, for all integers~$l$ and for
all $\om'\in \Om_l$. Hence, since $\int$ is left faithful, $\om=0$.
Therefore, $\s$ is injective, as required. \qed\enddemo

If $\int$ is a linear functional on a differential calculus, its {\it left kernel} is defined to be the set of all forms $\om$ for which ${\int \om'\om=0}$, for all $\om'\in \Om$. Obviously, the left kernel is a left ideal of $\Om$.
If the intersection of the left kernel of $\int$ with $\sa$ is the zero space, we say $\int$ is {\it weakly faithful}. Obviously, $\int$ is left faithful if, and only if, its left kernel is the zero space; hence, $\int$ is weakly faithful if it is left faithful, as one would expect.

\proclaim{Theorem 1.3} Let $\int$ be a twisted graded trace on a differential calculus $\omd$ over a unital algebra~$\sa$.

(1) $\int$ is weakly faithful if, and only if, for each element $a\in\sa$ for which
${\int a\om =0}$, for all $\om\in \Om$, we have $a=0$.

(2) If $\int$ is weakly faithful, then $\int$ admits exactly one twist automorphism. \endproclaim

\demo{Proof} First, suppose that $\int$ is weakly faithful. Let $\s$ be any twist automorphism of $\int$ and suppose that $a\in \sa$ and that ${\int a\om=0}$, for all $\om\in \Om$. Then ${\int \om \s^{-1}(a)=0}$. Hence, by weak faithfulness of $\int$, ${\s^{-1}(a)=0}$ and therefore, $a=0$.
This shows the forward implication in Condition~(1) and the reverse implication is shown by similar reasoning.

To see Condition~(2) holds, let $\r$ and $\s$ be twist automorphisms for $\int$. Then, for all $a\in \sa$ and $\om\in \Om$,
${\int (\r(a)-\s(a))\om}={\int \r(a)\om-\int \s(a)\om}={\int \om a-\int \om a=0}$. Hence, $\r(a)=\s(a)$. Using the fact that $\r d = d \r$ and $\s d = d \s$, it now follows immediately that $\r=\s$. \qed \enddemo

Let $N$ be a non-negative integer. We say that a linear functional $\int$ on $\Om$ is
{\it $N$-dimensional} if $\int\om=0$, for all $k$-forms, where $k\ne N$.

Suppose now $\int'$ is an $N$-dimensional, weakly faithful, closed twisted graded trace on a differential calculus $(\hat{\Om},d)$ over $\sa$ and let $\hat{\s}$ denote the  twist automorphism of $\int'$. We are going to construct a new, $N$-dimensional, differential calculus $\omd$ from $(\hat{\Om},d,\int')$ and a new, $N$-dimensional, closed twisted graded trace $\int$ on $\Om$ that is left faithful.

The twisted tracial property of $\int'$ implies that, for each form $\om\in \hat{\Om}$, the condition ${\int \om'\om=0}$, for all $\om'\in \hat{\Om}$, is equivalent to the condition ${\int \om\om'=0}$, for all $\om'\in \hat{\Om}$. Hence, if $I$ is the left kernel of~$\bint$, it is not only a left ideal of $\hat{\Om}$, but is also a right ideal. We denote by $\Om$ the quotient algebra $\hat{\Om}/I$. It is trivially verified that $\hat{\Om}_n\subseteq
I$ for all $n>N$ and that if $\om\in I$, then its $k$-th component
$\om_k$ belongs to $I$ also. It follows that if $\Om_k$ denotes
the image of $\hat{\Om}_k$ in the quotient algebra $\Om$, then
$\Om={\Om_0\oplus\cdots \oplus\Om_N}$. Moreover, this makes $\Om$
into a graded algebra. Since $I\cap \sa=0$, because $\bint$ is
weakly faithful, we may, and we do, identify $\Om_0$ with $\sa$.

If $\om'\in \hat{\Om}_k$ and $\om\in I$, then $\int'\om'd\om=(-1)^{k+1}\int' (d\om')\om=0$. This implies that $d\om \in I$. Hence, $d(I)\subseteq I$ and therefore $d$ induces a linear map ${d\colon \Om\to \Om}$. It is immediate that $d$ is a graded derivation on $\Om$ and, indeed, that $\omd$ is an $N$-dimensional differential calculus over $\sa$.

Since $\int'$ clearly annihilates $I$, we get an induced linear map $\int$ on $\Om$. Also, it is clear that $\hat{\s}(I)\subseteq I$, so that $\hat{\s}$ induces an algebra automorphism $\s$ on $\Om$. It is now easily verified that $\int$ is an $N$-dimensional, closed twisted graded trace on $\Om$ with $\s$ as its twist automorphism.

We call $\omd$ the {\it differential calculus} associated to $(\hat{\Om},d,\int')$ and $\int$ the {\it canonical} twisted graded trace on $\Om$.  The significant gains resulting from this construction are that $\omd$ is finite-dimensional and that $\int$ is left faithful.

It is straightforward to verify that if one starts with an $N$-dimensional differential calculus $\omd$ over $\sa$, and with a left faithful, closed twisted graded trace~$\int$ on~$\Om$, then (up to isomorphism) one can obtain $\Om$, $d$ and $\int$ by the preceding quotient construction from an $N$-dimensional, weakly faithful, closed twisted graded trace $\int'$ on $\bomd$.

The question now arises as to how we can obtain twisted graded traces on $\bomd$.
We shall see these arise from twisted cyclic cocycles. We shall discuss these objects and explain their relationship with twisted graded traces in Section~2.

Suppose now that $\sa$ is a unital $*$-algebra. We shall say that $\omd$ is a {\it $*$-differential calculus} over $\sa$ if it is a differential calculus over $\sa$ and if $\Om$ is endowed with a conjugate-linear map, ${\Om\to\Om}$, ${\om\mapsto \om^*}$, extending the involution on $\sa$, having the following properties:

(1) $(\om^*)^*=\om$, for all $\om\in \Om$;

(2) $(\om_1\om_2)^*=(-1)^{kl}\om_2^*\om_1^*$, for all $\om_1\in \Om_k$ and $\om_2\in \Om_l$;

(3) $d(\om^*)=(d\om)^*$, for all $\om\in\Om$.

We shall call the map, ${\om\mapsto \om^*}$, the {\it graded involution} of $\Om$. Notice that there is at most one such graded involution.

A linear map, ${\int \colon \Om\to \bC}$, is {\it self-adjoint} if $\int \om^*=(\int\om)^-$, for all $\om \in \Om$.

The universal differential calculus $\bomd$ of a $*$-algebra $\sa$ is a $*$-differential calculus in a natural way.
Suppose now $\int'$ is an $N$-dimensional, weakly faithful, self-adjoint, closed twisted graded trace on $\bomd$. Let $I$ be its left kernel, $\omd$ the associated $N$-dimensional differential calculus and $\int$ the canonical twisted graded trace on~$\Om$. Then $I$ is self-adjoint---that is, if $\om\in I$, then $\om^*\in I$---and $\omd$ is a $*$-differential calculus over $\sa$, where $(\om+I)^*=\om^*+I$, for all $\om\in \bar\Om$. To see $I$ is self-adjoint, suppose that $\om$ is a $k$-form belonging to $I$. If $\om'$ is an $(N\!-\!k)$-form, then ${\int' \om'\om^*}=({-1)^{k(N-k)}(\int' \om (\om')^*)^-=0}$. Hence, $\om^*\in I$. This proves $I^*\subseteq I$. It now follows easily that the involution $(\om+I)^*=\om^*+I$ makes $\omd$ into a $*$-differential calculus over $\sa$. It is equally easy to see that $\int$ is self-adjoint.

\subheading{2. Twisted cyclic cocycles and differential calculi}

Suppose that $\sa$ is a unital algebra. For $n\ge 0$, let $\bold C^n(\sa)$ denote the set of all multilinear maps from $\sa^{n+1}$ to $\bC$.
Set $\bold C^*(\sa)=\oplus_{n\in \bN} \bold C^n(\sa)$. Then $\bold C^*(\sa)$ is a graded linear space.
There exists a unique linear map, ${\bold b\colon \bold C^*(\sa)\to \bold C^*(\sa)}$, making $(\bold C^*(\sa),\bold b)$ a cochain complex for which, for $\f\in \bold C^n(\sa)$,
$$\align
(\bold b\f)(a_0,\dots,a_{n+1})=\sum_{i=0}^n (-1)^i\f(a_0,\dots,a_{i-1}&,a_ia_{i+1},a_{i+2},\dots,a_{n+1})\\
&+(-1)^{n+1}\f(a_{n+1}a_0,a_1,\dots,a_n).\endalign$$
The {\it Hochschild cohomology} $\bold H\bold H^*(\sa)$ of~$\sa$ is defined to be the cohomology of $(\bold C^*(\sa),\bold b)$. Thus, $\bold H\bold H^n(\sa)=\bold H^n(\bold C^*(\sa),\bold b)$ for all $n\in \bZ$.

The {\it permutation operator} $\l$ on $\bold C^*(\sa)$ is the linear isomorphism of degree zero, defined by setting $\l(\f)(a_0,a_1,\dots,a_n)=
(-1)^n\f(a_n,a_0,a_1,\dots,a_{n-1})$, for $n\ge 0$, $\f\in \bold C^n(\sa)$ and $a_0,\dots,a_n\in \sa$.
Set $\bold C^*_{\l}(\sa)=\oplus_{n\in \bN} \bold C^n_\l(\sa)$, where $\bold C^n_{\l}(\sa)={\{\f\in \bold C^n(\sa)\mid {\l}(\f)=\f\}}$. The coboundary operator $\bold b$ leaves  each space $\bold C^n_{\l}(\sa)$ invariant and therefore its restriction makes $(\bold C^*_{\l}(\sa),\bold b)$
into a cochain complex. The cohomology of this complex is denoted by $\bold H^*_{\l}(\sa)$ and called the {\it cyclic cohomology} of $\sa$. Thus, $\bold H^n_{\l}(\sa)=\bold H^n(\bold C^*_{\l}(\sa),\bold b)$.

It will be useful to recall also the degree 1 operator $\bold b'$
on $\bold C^*(\sa)$ defined by the formula $$(\bold
b'\f)(a_0,\dots,a_{n+1})=\sum_{i=0}^n
(-1)^i\f(a_0,\dots,a_{i-1},a_ia_{i+1},a_{i+2},\dots,a_{n+1})$$ for $n \geq 0$ and $\f \in \bold C^n(\sa)$. It
is well known that $(\bold b')^2=0$ and that the cohomology of the
cochain complex ${(\bold C^*(\sa),\bold b')}$ is trivial, ${\bold H^*(\bold
C^*(\sa),\bold b')=0}$.

We generalize the definition of cyclic cohomology now. Suppose
that $(\sa,\s)$ is a pair consisting of a unital algebra $\sa$ and an
algebra automorphism ${\s\colon \sa\to \sa}$. We get a new operator
corresponding to the permutation operator, a linear isomorphism
${{\l}\colon \bold C^*(\sa)\to \bold C^*(\sa)}$ of degree zero, by
setting $${\l}(\f)(a_0,a_1,\dots,a_n)=
(-1)^n\f(\s(a_n),a_0,a_1,\dots,a_{n-1})$$ for $n \geq 0$ and $\f \in \bold C^n(\sa)$. We set $\bold
C^*_{\l}(\sa,\s)=\oplus_{n\in \bN} \bold C^n_{\l}(\sa,\s)$, where
$\bold C^n_{\l}(\sa,\s)= {\{\f\in \bold C^n(\sa)\mid {\l}(\f)=\f\}}$.
We shall make $\bold C^*_{\l}(\sa,\s)$ into a cochain complex whose
cohomology will be a ``twisted'' version of ordinary cyclic
cohomology. To this end we introduce new operators $\bold c$ and
$\bold b$ on $\bold C^*(\sa)$, both of degree 1. These are defined
by setting $\bold b=\bold b'+\bold c$, where, for $\f\in \bold
C^n(\sa)$, and ${a_0,\dots,a_n\in \sa}$, $$(\bold
c\f)(a_0,\dots,a_{n+1})=(-1)^{n+1}\f(\s(a_{n+1})a_0,a_1,\dots,a_n).$$
Thus, $\bold b$ is a ``twisted'' version of the usual Hochschild
coboundary operator. To see that $\bold b^2=0$, one uses the fact
that $(\bold b')^2=0$ and proves the easily verified fact that
${\bold c \bold b'+\bold b'\bold c+\bold c^2=0}$. As in the
classical cyclic cocycle theory, one can show that $\bold
b'(1-{\l})=(1-{\l})\bold b$. This immediately implies that $\bold
C^*_{\l}(\sa,\s)= {\{\f\in \bold C^*(\sa)\mid {\l}\f=\f\}}$ is
invariant under $\bold b$. Hence, by restricting $\bold b$, we get
a cochain complex $(\bold C^*_{\l}(\sa,\s),\bold b)$. We denote by
$\bold H^*_{\l}(\sa,\s)$ the cohomology of this complex and call it the
{\it twisted cyclic cohomology} of $(\sa,\s)$. We denote by
$\bold Z^n_{\l}(\sa,\s)$ and $\bold B^n_{\l}(\sa,\s)$ the $n$-cocyles and
$n$-coboundaries for the complex $(\bold C^*_{\l}(\sa,\s),\bold b)$.
We call the elements of these spaces the {\it twisted cyclic
$n$-cocyles} and {\it $n$-coboundaries} of $(\sa,\s)$, respectively.

Clearly, if $\s=\id_\sa$, then $\bold H^*_{\l}(\sa,\s)=\bold H^*_{\l}(\sa)$.

\proclaim{Theorem 2.1} Let $\omd$ be a differential calculus over a unital algebra~$\sa$ and suppose that $\int$ is an $N$-dimensional, closed, twisted graded trace on $\Om$. Define the function, ${\f\colon \sa^{N+1}\to \bC}$, by setting
$${\f(a_0,\dots,a_N)}={\int a_0da_1\cdots da_N}.$$ Let $\s$ be an automorphism of~$\sa$ for which ${\int \s(a)\om}={\int \om a}$, for all $a\in \sa$ and $\om\in \Om_N$. Then $\f$ belongs to $\bold Z^N_{\l}(\sa,\s)$. \endproclaim

\demo{Proof} We show first that ${\l} \f=\f$. Let ${a_0,\dots a_N}$ be elements of $\sa$. Then, since $\int$ is closed, and ${da_0\cdots da_{N-1}}={d(a_0da_1\cdots da_{N-1})}$, we have
$$\align {\l} \f(a_0,\dots,a_N)
&=(-1)^N\int \s(a_N)da_0\cdots da_{N-1}
=(-1)^N \int (da_0\cdots da_{N-1})a_N
\\&=\int a_0(da_1\cdots da_{N-1})da_N
=\f(a_0,\dots,a_N).\endalign$$ To show that $\bold b\f=0$, we
shall use the fact that $$\sum_{i=1}^N (-1)^i da_1\cdots
d(a_ia_{i+1})\cdots da_{N+1}= (-1)^N (da_1\cdots da_N)a_{N+1}-
a_1da_2\cdots da_{N+1}, \quad (2.1) $$ for all $a_1,\dots,a_{N+1}\in \sa$ (this
is well known, see~\cite{\CON, p. 187}). It follows from this
equality, and from the twisted tracial property of $\int$, that
$$\align \bold b\f(a_0,\dots,a_{N+1}) & =\sum_{i=1}^N (-1)^i\int
a_0da_1\cdots d(a_ia_{i+1})\cdots da_{N+1}
\\&\phantom{xxxxxxx}+\int a_0a_1da_2\cdots da_{N+1}
+(-1)^{N+1}\int \s(a_{N+1})a_0da_1\cdots da_N\\ &=\int
a_0((-1)^N(da_1\cdots da_N)a_{N+1}-a_1da_2\cdots da_{N+1})
\\&\phantom{xxxxxxx}+\int a_0a_1da_2\cdots da_{N+1}+(-1)^{N+1}\int a_0(da_1\cdots da_N)a_{N+1}=0.\endalign$$
The theorem is now proved \qed\enddemo

We call $\f$ the twisted cyclic cocycle {\it associated to} $\omd$ and $\int$.

\proclaim{Theorem 2.2} Let $\s$ be an automorphism of a unital algebra~$\sa$ and let $\f\in \bold Z^N_{\l}(\sa,\s)$, for some integer $N\ge 0$. Then there exists an $N$-dimensional differential calculus $\omd$ over~$\sa$ and an
$N$-dimensional, closed twisted graded trace $\int$ on $\Om$ such that $\f$ is the twisted cyclic cocycle associated to $\omd$ and $\int$. \endproclaim

\demo{Proof} Define an $N$-dimensional linear functional $\bint$ on the universal differential calculus~$\bom$ over~$\sa$ by setting ${\bint a_0da_1\cdots da_N}=\f(a_0,\dots,a_N)$ and ${\bint da_1\cdots da_N}=  0$, for all ${a_0,\dots,a_N\in \sa}$. By definition, $\bint$ is closed.

Next we show that ${\bint \om a_{N+1}}={\bint \s(a_{N+1}) \om}$, for all $a_{N+1}\in \sa$ and all $\om\in \bom$. Clearly, to show this, we may suppose that $\om=a_0da_1\cdots da_N$ or $\om = da_1\cdots da_N$, for some elements ${a_0,\dots,a_N\in \sa}$. Then, using the fact that $\bold b\f=0$ and therefore, $\bold b'\f =-\bold c\f$, and again using Equation~(2.1), we have
$$\align
\bint & \s(a_{N+1})a_0da_1\cdots da_N
= (-1)^{N+1}\bold c\f(a_0,\dots,a_{N+1})
= (-1)^N \bold b'\f(a_0,\dots,a_{N+1})
\\&= (-1)^N  \sum_{i=0}^N (-1)^i\f(a_0,\dots,a_ia_{i+1},\dots,a_{N+1})
\\&= (-1)^N (\sum_{i=1}^N (-1)^i \bint a_0da_1\cdots d(a_ia_{i+1})\cdots da_{N+1}
+ \bint a_0a_1da_2\cdots da_{N+1})
\\&= (-1)^N (\bint \! a_0((-1)^N (da_1\cdots da_N)a_{N+1}-a_1da_2\cdots da_{N+1})
 +\bint \! a_0a_1da_2\cdots da_{N+1})
\\&=  \bint a_0(da_1\cdots da_N)a_{N+1}.
\endalign$$
In the other case
$$\align  \bint & \s(a_{N+1})da_1\cdots da_N = \f(\s(a_{N+1}),a_1,\dots,a_N)
= (-1)^N \f(a_1,\dots,a_N,a_{N+1}) \\
& = (-1)^N \bint a_1 da_2 \cdots da_{N+1} = \bint (da_1 \cdots da_N) a_{N+1} \ ,
\endalign$$
where we used the closedness of $\bint$ and the aforementioned fact in the last equality.

It follows now that $\bint$ is a twisted graded trace. Now let
$\omd$ be the $N$-dimensional differential calculus obtained from
$\bom$ by truncation, and let $\int$ be the restriction of $\bint$
to $\Om$. Clearly, $\int$ is again a closed twisted graded trace
and $\f$ is the twisted cyclic cocycle associated to $\omd$ and
$\int$. \qed\enddemo

If $\s$ is an automorphism of a unital algebra~$\sa$ and $\f\in \bold C^*_{\l}(A,\s)$, we say that $\f$ is {\it left faithful} if, for each element $a$ in $\sa$, we have $a=0$, if ${\f(aa_0,a_1,\dots a_N)=0}$, for all ${a_0,\dots,a_N\in \sa}$. Since $\l \f=\f$, we have, for each index ${i=0,\dots,N}$, $a=0$, if
${\f(a_0,a_1,\dots,aa_i,\dots a_N)=0}$, for all ${a_0,\dots,a_N\in \sa}$.

\proclaim{Theorem 2.3} Let $\s$ be an automorphism of a unital algebra~$\sa$ and let $\f\in \bold Z^N_{\l}(\sa,\s)$, for some integer $N\ge 0$. If $\f$ is left faithful, then there exists an
$N$-dimensional  differential calculus $\omd$ over~$\sa$ and a left faithful
$N$-dimensional, closed twisted graded trace $\int$ on $\Om$ such that $\f$ is the twisted cyclic cocycle associated to $\omd$ and $\int$. \endproclaim

\demo{Proof} Define an $N$-dimensional linear functional $\bint$ on the universal differential calculus~$\bom$ over~$\sa$ by setting ${\bint a_0da_1\cdots da_N}=\f(a_0,\dots,a_N)$ and ${\bint da_1\cdots da_N}=  0$, for all ${a_0,\dots,a_N\in \sa}$. We saw in the proof of the preceding theorem that $\bint$ is a closed twisted graded trace. The faithfulness assumption on $\f$ ensures that $\bint$ is weakly faithful.
Now let $\omd$ be the $N$-dimensional calculus associated to $\bom$ and $\bint$ and let $\int$ be the canonical $N$-dimensional, left faithful, closed twisted graded trace on $\Om$. Then $\f$ is clearly the twisted cyclic cocycle associated to~$\int$. \qed\enddemo

To round off this circle of ideas, let us note that if $\int$ is any $N$-dimensional, weakly faithful, closed twisted graded trace on a differential calculus~$\omd$ over a unital algebra~$\sa$, the associated twisted cyclic cocycle~$\f$ is clearly left faithful.

We turn now to the case of $*$-differential calculi. If $\omd$ is such a calculus over a unital $*$-algebra~$\sa$, then it is readily verified that, for all 1-forms ${\om_1,\dots,\om_N}$ of $\Om$, ${(\om_1\cdots\om_N)^*}={s_N \om^*_N\cdots\om^*_1}$, where $(s_N)$ is the sequence of scalars defined inductively by $s_1=1$ and $s_{N+1}=(-1)^Ns_N$.
If $\f$ is the $N$-cocycle associated to an $N$-dimensional weakly faithful, closed, self-adjoint, twisted graded trace $\int$ on $\Om$, then $\f^*=\f$, where ${\f^*(a_0,\dots,a_N)}=
{s_{N+1}\bar \f(a_N^*,\dots,a^*_0)}$ (as usual, $\bar\f$ is the complex congugate function corresponding to $\f$, so that $\bar\f(x)=\overline{\f(x)}$).
To see that $\f^*=\f$, observe that, if $\s$ is a twist automorphism associated to~$\int$, then
$$\align
\f^*&(a_0,\dots,a_N)=
{s_{N+1}\bar\f(a_N^*,\dots,a^*_0)}={s_{N+1}(-1)^N\bar\f(\s(a_0^*),a_N^*,\dots,a_1^*)}\\
&={s_{N+1}(-1)^N(\int \s(a_0^*)(da_N^*)\cdots(da_1^*))^-}
={(-1)^Ns_Ns_{N+1} \int (da_1)\cdots (da_n)\s(a_0^*)^*}\\
&={s^2_{N+1}\int (da_1)\cdots (da_N)\s^{-1}(a_0)}={\int a_0da_1\cdots da_N}=
{\f(a_0,\dots,a_N)}.
\endalign$$
Here, in the third last equation, we have used the easily verified fact that $\s^{-1}(a^*)=\s(a)^*$, for all $a\in \sa$ (this uses weak faithfulness of $\int$).

These observations motivate the following definitions.

If the function, ${\f\colon\sa^{N+1}\to \bC}$, is multilinear, we define $\f^*$ by setting ${\f^*(a_0,\dots,a_N)}={s_{N+1}\bar \f(a_N^*,\dots,a^*_0)}$, for all ${a_0,\dots,a_N\in \sa}$.

If $\s$ is an automorphism of $\sa$ such that $\s(a)^*=\s^{-1}(a^*)$, for all $a\in\sa$, then we call $\s$ {\it regular}. As we observed above,
the restriction to $\sa$ of a twist automorphism associated to a weakly faithful, self-adjoint twisted graded trace is regular. Another observation: if $\s$ is any self-adjoint automorphism of~$\sa$ and $\s^2=\id$, then $\s$ is regular.

It is easy check that, if $\s$ is any regular automorphism of $\sa$, and ${\f\in \bold C^N_\l(\sa,\s)}$, then ${\f^*\in \bold C^N_\l(\sa,\s)}$. It is also the case that, if $\bold b\f=0$, then $\bold b\f^*=0$. However, this requires some proof, so we give the details. It clearly sufices to show that, if ${a_0,\dots,a_{N+1}\in \sa}$,
then
$$\sum_{i=0}^N (-1)^i \f(a_{N+1}^*,\dots,a_{i+1}^*a_i^*,\dots,a_0^*)+(-1)^{N+1}\f(a_N^*,\dots,a_1^*,a_0^*\s(a_{N+1})^*)=0.$$
Set $b_i=a_{N+1-i}^*$, for ${i=0,\dots,N+1}$. Multiplying the above equation
by $(-1)^N$ and using the fact that $\s(a_{N+1})^*=\s^{-1}(a_{N+1}^*)=\s^{-1}(b_0)$,
we see that we need only show that
$$\sum_{i=0}^N (-1)^{N-i} \f(b_0,\dots,b_{N-i}b_{N-i+1},\dots,b_{N+1})+
(-1)^{2N+1}\f(b_1,\dots,b_N,b_{N+1}\s^{-1}(b_0))=0.$$
Now we use the fact that $\l\f=\f$, which implies that
${(-1)^N\f(b_1,\dots,b_N,b_{N+1}\s^{-1}(b_0))}={\f(\s(b_{N+1})b_0,b_1,\dots,b_N)}$, to see that we have only to show that
$$\sum_{i=0}^N (-1)^{N-i} \f(b_0,\dots,b_{N-i}b_{N-i+1},\dots,b_{N+1})+
(-1)^{N+1}{\f(\s(b_{N+1})b_0,b_1,\dots,b_N)}=0;$$
that is, it suffices to show that
$$\sum_{i=0}^N (-1)^i \f(b_0,\dots,b_ib_{i+1},\dots,b_{N+1})+
(-1)^{N+1}{\f(\s(b_{N+1})b_0,b_1,\dots,b_N)}=0.$$
However, this is true, since it is just the equation
${(\bold b'+\bold c)\f(b_0,\dots,b_{N+1})=0}$, which holds because $\bold b\f=0$, by assumption.

If we define $\f$ to be {\it self-adjoint}, if $\f^*=\f$, then the preceding observations, together with the easily  checked equation $(\f^*)^*=\f$, show that every element ${\f\in \bold Z^N_\l(\sa,\s)}$ can be written in the form $\f=\f_1+i\f_2$, for some self-adjoint elements $\f_1$ and $\f_2$ in ${\bold Z^N_\l(\sa,\s)}$. (Of course, one sets ${\f_1=(\f+\f^*)/2}$ and ${\f_2=(\f-\f^*)/2i}$.)

Now suppose that $\int$ is an $N$-dimensional, closed, twisted graded trace on a $*$-differential calculus $\omd$. If the twisted cyclic $N$-cocycle $\f$ associated to $\int$ is self-adjoint, then $\int$ is self-adjoint. To see this we need only show that ${(\int \om)^-}={\int \om^*}$, where ${\om =a_0da_1\cdots da_N}$ or ${\om = da_1\cdots da_N}$, for elements ${a_0,\dots,a_N}$ belonging to $\sa$. However, we have
$$\align
(&\int\om)^-=\bar\f(a_0,\dots,a_N)=\bar\f^*(a_0,\dots,a_N)
=s_{N+1}\f(a_N^*,\dots,a_0^*)\\
&=s_{N+1}\int a_N^*(da_1^*)\cdots (da_0^*)
=(-1)^N\int ((da_0)\cdots (da_{N-1})a_N)^*\\
&=(-1)^N\int (d(a_0da_1\cdots da_{N-1})a_N)^*
=\int (a_0d(a_1da_2\cdots da_N))^*=\int \om^*.
\endalign$$
In the second last equation we used the fact that ${\int d=0}$ and that
${d((a_0da_1\cdots da_{N-1})a_N)}\allowmathbreak={d(a_0da_1\cdots da_{N-1})a_N+(-1)^{N-1}a_0d(a_1da_2\cdots da_N)}$.

If $\om = da_1\cdots da_N$, it is clear that $\int \om = 0 = \int \om^*$ due to the closedness of $\int$.

We sum up our observations in the following theorem.

\proclaim{Theorem 2.4} Let $\sa$ be a unital $*$-algebra and let $\s$ be a regular (algebra) automorphism of $\sa$. Let $\int$ be an $N$-dimensional, closed, twisted graded trace
on a $*$-differential calculus $\omd$ over~$\sa$, and suppose that its twist automorphism extends $\s$. Let $\f$ be the twisted cyclic $N$-cocycle associated to $\int$, so that ${\f\in \bold Z^N_\l(\sa,\s)}$. Then $\f$ is self-adjoint if, and only if, $\int$ is self-adjoint. \endproclaim

\subheading{3. Twisted cyclic cohomology}

In this section we briefly consider the twisted cyclic cohomology theory of a pair $(\sa,\s)$, where $\sa$ is a unital algebra and $\s$ is an automorphism of $\sa$. We shall be particularly interested in the construction of analogues of the important operators $\bold S$ and $\bold B$ occurring in the classical cyclic cohomology theory. These are used to relate twisted cyclic cohomology to twisted Hochschild cohomology. We begin by defining the latter. Note that if $\f\in \bold C^n(\sa)$, then $(\l^{n+1}\f)(a_0,\dots,a_n)
=\f(\s(a_0),\dots,\s(a_n))$, for all $a_0,\dots,a_n\in \sa$.
Let $\bold C^*(\sa,\s)= \oplus_{n\in \bN} \bold C^n(\sa,\s)$, where $\bold C^n(\sa,\s)={\{\f\in \bold C^n(\sa)\mid \l^{n+1}\f = \f \}}$.
One can show that, for $\f\in \bold C^n(\sa)$, we have $\bold b\l^{n+1}\f=\l^{n+2}\bold b\f$ and $\bold b'\l^{n+1}\f=\l^{n+2}\bold b'\f$. It follows that $\bold C^*(\sa,\s)$ is invariant for $\bold b$ and $\bold b'$ and therefore we get a cochain complex $(\bold C^*(\sa,\s),\bold b)$. We denote its cohomology by $\bold H\bold H(\sa,\s)$ and call it the {\it twisted Hochschild cohomology} of the pair~$(\sa,\s)$.

We shall now get the twisted cyclic cohomology as the cohomology of the total complex of a bicomplex. To define this bicomplex we introduce the operator $\bold N$ of degree zero on $\bold C^*(\sa,\s)$, defined, for $\f\in \bold C^n(\sa,\s)$, by setting $\bold N\f=\sum_{i=0}^n\l^i\f$. One can show that $\bold b \bold N=\bold N\bold b'$ and
$(1-\l)\bold b=\bold b'(1-\l)$ and $\bold N(1-\l)=0$. Hence, for $\bold C^n=\bold C^n(\sa,\s)$, the following diagram defines a bicomplex
$$\matrix
\vdots &  & \ \vdots &  &  \vdots &  & \ \vdots &  & \\
\bold b\uparrow \quad &  & -\bold b'\uparrow \quad &  & \bold b\uparrow \quad &  & -\bold b'\uparrow \quad &  & \\
\bold C^2 & \mapp{1-\l} & \bold C^2 & \mapp{\bold N} & \bold C^2 & \mapp{1-\l} & \bold C^2 & \mapp{\bold N} &\cdots\\
\bold b\uparrow \quad &  & -\bold b'\uparrow \quad &  & \bold b\uparrow \quad &  & -\bold b'\uparrow \quad &  & \\
\bold C^1 & \mapp{1-\l} & \bold C^1 & \mapp{\bold N} & \bold C^1 & \mapp{1-\l} & \bold C^1 & \mapp{\bold N} &\cdots\\
\bold b\uparrow \quad &  & -\bold b'\uparrow \quad &  & \bold b\uparrow \quad &  & -\bold b'\uparrow \quad &  & \\
\bold C^0 & \mapp{1-\l} & \bold C^0 & \mapp{\bold N} & \bold C^0 & \mapp{1-\l} & \bold C^0 & \mapp{\bold N} &\cdots\\
\endmatrix$$
We denote this bicomplex by $\bold C^{**}(\sa,\s)$ and its total complex by $\bold T^*(\sa,\s)$. The entry in the bicomplex at the position $(m,n)$ is $\bold C^{m,n}(\sa,\s)=\bold C^n(\sa,\s)$. We denote the cohomology of $\bold T^*(\sa,\s)$ by $\bold H\bold C^*(\sa,\s)$. We shall see that this is isomorphic to $\bold H_\l(\sa,\s)$. The advantage of this alternative description is that it enables us to define the operators $\bold S$ and $\bold B$ in a natural way.

We define a cochain map $\pi$ from the complex $\bold C^*_\l(\sa,\s)$ to the complex $\bold T^*(\sa,\s)$ by mapping $x$ in $\bold C^n_\l(\sa,\s)$ onto $(x,0,\dots,0)$ in $\bold T^n(\sa,\s)=\oplus_{i=0}^n \bold C^{i,n-i}(\sa,\s)$. Then one can show that the induced linear map, ${\pi_*\colon \bold H^*_\l(\sa,\s)\to \bold H\bold C^*(\sa,\s)}$, is an isomorphism.

We now define $\bold C^{**}_{[2]}$ to be the cochain bicomplex obtained from $\bold C^{**}(\sa,\s)$ by restricting to the first two columns and setting all other columns equal to zero. Let $\bold T^*_{[2]}(\sa,\s)$ be the total complex of $\bold C^{**}_{[2]}$. We define a cochain map ${\ta}$ from $\bold T^*_{[2]}(\sa,\s)$ to $\bold C^*(\sa,\s)$ by setting ${\ta}(x)=x$, for $x$ in $\bold T^0_{[2]}(\sa,\s)=\bold C^0(\sa,\s)$ and setting ${\ta}(x_0,x_1)=x_0$, for $(x_0,x_1)$ in $\bold T^n_{[2]}(\sa,\s)=\bold C^n(\sa,\s)\oplus \bold C^{n-1}(\sa,\s)$, where $n>0$. The induced map, ${{\ta}_*}$ mapping ${\bold H^*(\bold T^*_{[2]}(\sa,\s))}$ to ${\bold H\bold H^*(\sa,\s)}$, is an isomorphism.

Now we define a cochain map of degree 2 on $\bold T^*(\sa,\s)$ by shifting its chain bicomplex two columns to the right; more precisely, if $x=(x_0,\dots,x_n)\in \bold T^n(\sa,\s)$, set ${\bold R}(x)=(0,0,x_0,\dots,x_n)$. Let ${\bold P}$ be the degree zero cochain map from $\bold T^*(\sa,\s)$ to $\bold T^*_{[2]}(\sa,\s)$ obtained by projecting; more precisely, ${\bold P}(x)=x$ for $x\in \bold T^0(\sa,\s)$ and ${\bold P}(x)=(x_0,x_1)$, for $x=(x_0,\dots,x_n)\in \bold T^n(\sa,\s)$, where $n>0$. This gives a short exact sequence of cochain maps
$$0\to \bold T^*(\sa,\s)\mapp{\bold R} \bold T^*(\sa,\s)\mapp{\bold P} \bold T^*_{[2]}(\sa,\s)\to 0.$$
On the cohomological level we therefore get an exact triangle
$$\triangle{\bold H^*(\bold T^*_{[2]}(\sa,\s))}\partial{\bold H^*(\bold T^*(\sa,\s))}{{\bold R}_*}{\bold H^*(\bold T^*(\sa,\s))}{{\bold P}_*}$$
Finally, we define the linear maps ${\bold I\colon \bold H^*_\l(\sa,\s)\to \bold H\bold H^*(\sa,\s)}$,
${\bold S\colon \bold H^*_\l(\sa,\s)\to \bold H^*_\l(\sa,\s)}$ and
${\bold B\colon \bold H\bold H^*(\sa,\s)\to \bold H_\l^*(\sa,\s)}$ of degrees 0, 2 and -1 respectively by setting $\bold I={\ta}_*{\bold P}_*\pi_*$, $\bold S=\pi_*^{-1}{\bold R}_*\pi_*$ and
$\bold B=\pi_*^{-1}\partial{\ta}_*^{-1}$. This gives us an exact triangle
$$\triangle{\bold H\bold H^*(\sa,\s)}{B}{\bold H^*_\l(\sa,\s)}{\bold S}{\bold H^*_\l(\sa,\s)}{\bold I}$$
By expansion of this we get a long exact sequence
$$\dots\to \bold H^{n-2}_\l(\sa,\s)\mapp{\bold S} \bold H^n_\l(\sa,\s)\mapp{\bold I} \bold H\bold H^n(\sa,\s)
\mapp{\bold B} \bold H^{n-1}_\l(\sa,\s) \mapp{\bold S} \bold H^{n+1}_\l(\sa,\s)\to\cdots$$
Thus, we have indicated how the principal results of the elementary theory of cyclic cohomology extends to the twisted case. Since the proofs in this more general setting are essentially the same as in the non-twisted case, we have omitted the details.

\subheading{4. Left-covariant differential calculi}

Differential calculi that are left-covariant are of prime importance for the theory. We shall introduce this concept now.
For this we need to suppose that $\sa$ is endowed with a co-multiplication $\D$ making the pair $\sad$ a Hopf algebra (such an algebra is unital by assumption). In the sequel we shall use a number of elementary results about Hopf algebras without explicit reference. A good general source for this material is~\cite{\ABE}.

Recall that a left-covariant bi-module over $\sad$ is a pair $(\Gamma,\D_\Gamma)$, where $\Gamma$ is a bi-module over $\sa$, and $\D_\Gamma$ is a linear map from $\Gamma$ to${\sa\otimes\Gamma}$ such that the following conditions hold:

(1) ${(\D\otimes \id_\Gamma)\D_\Gamma}={(\id_\sa\otimes \D_\Gamma)\D_\Gamma}$ and ${(e\otimes \id_\Gamma)\D_\Gamma}=\id_\Gamma$, where $e$ is the co-unit of $(\sa,\D)$, (that is, $\D_\Gamma$ is a left co-action);

(2) $\D_\Gamma(a\g b)=\D(a)\D_\Gamma(\g)\D(b)$, for all $\g\in \Gamma$ and $a,b\in \sa$.

An element $\g\in \Gamma$ is said to be {\it left invariant} if $\D_\Gamma(\g)={1\otimes \g}$. We denote by $\Gamma^{\text{inv}}$ the linear space of left-invariant elements of $\Gamma$.

If $a\in \sa$ and $f$ is a linear functional on $\sa$, we set $f*a={(\id_\sa\otimes f)\D(a)}$.
We shall make use of the following result from the theory of left-covariant bi-modules.

\proclaim{Theorem 4.1 (S.L. Woronowicz \cite{\WORA, \WORC})} Let $(\Gamma,\D_\Gamma)$ be a left-covariant bi-module over a Hopf-algebra $(\sa,\D)$.

(1) There is a unique isomorphism of left $\sa$-modules from ${\sa\otimes \Gamma^{\text{inv}}}$ onto $\Gamma$ that maps ${a\otimes \g}$ onto $a\g$, for all $a\in \sa$ and $\g\in \Gamma^{\text{inv}}$.

(2) Suppose that the family of elements $(\g_i)_{i\in I}$ is a linear basis for $\Gamma^{\text{inv}}$. Then it is a free left $\sa$-module basis for $\Gamma$ and also a free right $\sa$-module basis of $\Gamma$. Moreover, there exist linear functionals $f_{jk}$ on $\sa$, for all $j,k\in I$, such that $f_{jk}(ab)=\sum_{i\in I} f_{ji}(a)f_{ik}(b)$ and $f_{jk}(1)=\d_{jk}$ and for which we have the equations
$\g_ja={\sum_{i\in I}(f_{ji}*a)\g_i}$ and $a\g_j=
{\sum_{i\in I} \g_i((f_{ji}\k^{-1})*a)}$, where $\k$ is the co-inverse for $\sad$.\endproclaim

When we consider a sum $\sum_{i\in I} x_i$ of a family $(x_i)_{i\in I}$ of elements in a vector space $X$ with no topological structure, it is understood that $x_i=0$ for all but a finite number of indices $i\in I$.

Let $\omd$ be a unital differential calculus over $\sa$ such that $d1 = 0$. This is a bi-module over $\sa$ in a natural way. If the map, ${\dom\colon\Om\to \sa\otimes \Om}$, makes $\Om$ into a left-covariant bi-module and
${(\id_\sa\otimes d)\dom}=\dom d$, and $\dom(a)=\D(a)$, for all $a\in \sa$, we call the triple ${(\Om,d,\dom)}$ a {\it left-covariant differential calculus} over $\sad$. A moment's reflection, using the fact that $\Om$ is generated as an algebra by the elements $a$ and $da$, where $a\in \sa$, shows that only one such left action $\dom$ can exist making ${(\Om,d,\dom)}$ a left-covariant calculus. For this reason, we often speak of the left-covariant differential calculus $\omd$, omitting explicit reference to $\dom$. Henceforth, we shall also often speak of the Hopf algebra~$\sa$, omitting explicit reference of the co-multiplication $\D$.

The map $\dom$ is automatically of degree zero, where we regard ${\sa\otimes \Om}$ as graded algebra
in the obvious way (its space of $k$-forms is the tensor product ${\sa\otimes \Om_k}$).

The linear span of the set ${\D(\sa)(\sa\otimes 1)}=
{\{\D(a)(b\otimes 1)}\mid {a,b\in \sa\}}$ is equal to ${\sa\otimes \sa}$ (this is true for any Hopf algebra). It follows from this that the linear span of ${\dom(\Om)(\sa\otimes 1)}$ is equal to ${\sa\otimes \Om}$.

We shall denote the linear space of left-invariant $k$-forms of $\Om$ by $\ominv_k$.

\smallskip

Let $\sa$ be any unital algebra (not necessarily the underlying algebra of a Hopf algebra). In section 1 we introduced the universal differential algebra $\bomd$ over $\sa$ (which is not unital). But there also exists a universal unital differential algebra over $\sa$ and this is the one we will be working with in the rest of this paper. There exists a unital graded differential algebra $\tomd$, for which
$\tom_0=\sa$, that has the following universal property: If
$\s$ is a unital algebra homomorphism from $\sa$ into the algebra
$\Om_0$ of $0$-forms of a unital graded differential algebra $\omd$, then
there exists a unique unital algebra homomorphism $\tilde{\s}$ from $\tom$ to
$\Om$ extending $\s$ such that $\tilde{\s} d=d\tilde{\s}$. This property
uniquely determines $\tomd$ (up to isomorphism). Note that $d1=0$.

We shall use the following useful property of $\tomd$:

Let $n$ be a non-negative integer and $T$ a multilinear map from ${\sa^{n+1}}$ to a linear space $Y$ such that ${T(a_0,\dots,a_n)=0}$, if any of the elements ${a_1,\dots,a_n}$ is a scalar. Then there is a unique linear map $\hat T$ from $\tom_n$ to $Y$ for which $\hat T(a_0da_1\cdots da_n)=T(a_0,a_1,\dots,a_n)$, for all $a_0,a_1,\dots,a_n\in \sa$.

Theorem 1.1 remains valid for $\tomd$ in place of $\bomd$, provided $\s_0$ is assumed to be unital.

If $(\sa,\D)$ is a Hopf algebra, then the universal unital calculus $\tomd$ over $\sa$ is a left-covariant calculus over $\sad$. To see this, first observe that ${\sa\otimes \tom}$ can be made into a differential calculus, where ${\id_\sa\otimes d}$ is its differential. The map $\D$, regarded as an algebra homomorphism
from $\sa$ to the 0-forms of ${\sa\otimes \tom}$, extends to an algebra homomorphism $\D'$ from
$\tom$ to ${\sa\otimes \tom}$ such that $\D'd={(\id_\sa\otimes d)\D'}$. It now follows from the next lemma that ${(\tom,d,\D')}$ is a left-covariant differential calculus over $\sad$.

\proclaim{Lemma 4.2} Let $\omd$ be a unital differential calculus over a Hopf algebra $\sad$ such that $d 1 = 0$ and suppose that ${\D_\Om\colon \Om\to \sa\otimes \Om}$ is an algebra homomorphism extending ${\D\colon \sa\to \sa\otimes \sa}$ such that ${(\id_\sa\otimes d)\D_\Om}=\D_\Om d$. Then ${(\Om,d,\D_\Om)}$ is a left-covariant differential calculus. \endproclaim

\demo{Proof} We have to prove that ${(\D\otimes \id_\Om)\D_\Om}={(\id_{\sa}\otimes \D_\Om)\D_\Om}$ and ${(e\otimes \id_\Om)\D_\Om}=\id_\Om$, where $e$ is the co-unit of $(\sa,\D)$.
We shall prove only the first of these equations; the proof of the second is straightforward. Since
${(\D\otimes \id_\Om)\D_\Om}$ and ${(\id_{\sa}\otimes \D_\Om)\D_\Om}$ are homomorphisms and $\Om$ is generated as an algebra by the forms $a$ and $da$, where $a\in \sa$, we need only see that these homomorphisms are equal at such forms. This is obvious in the case of the elements $a$, since $\D_\Om(a)=\D(a)$. For $da$ we have
$$\align
(\D&\otimes \id_\Om)\D_\Om d(a)=
{(\D\otimes \id_\Om)(\id_{\sa}\otimes d)\D(a)}=
{(\id_{\sa}\otimes\id_{\sa}\otimes d)(\D\otimes\id_\sa)\D(a)}\\
&={(\id_{\sa}\otimes\id_{\sa}\otimes d)(\id_{\sa}\otimes\D)\D(a)}=
{(\id_{\sa}\otimes \D_\Om d)\D(a)}\\&=
{(\id_{\sa}\otimes \D_\Om)(\id_{\sa}\otimes d)\D(a)}
={(\id_{\sa}\otimes \D_\Om)\D_\Om d(a)}.
\endalign$$
This proves the lemma.\qed\enddemo

Recall that a linear functional $h$ on a Hopf algebra $\sa$ is said to be {\it left-invariant} if ${(\id\otimes h)\D(a)=h(a)1}$, for all $a\in \sa$, where 1 is the unit of $\sa$. Similarly, a linear functional $h'$ on $\sa$ is {\it right-invariant} if ${(h'\otimes \id)\D(a)=h'(a)1}$, for all $a\in \sa$.
Such functionals do not necessarily exist. It is easily seen that there is at most one unital linear functional $h$ on $\sa$ that is both left and right invariant. We call such a functional a {\it Haar integral} of $\sa$. In the sequel, we shall be principally interested in working with Hopf algebras that admit Haar integrals. If $\sa$ is the Hopf algebra associated to a compact quantum group in the sense of Woronowicz, then it admits a Haar integral. From the point of view of relevance of the theory we are developing here, the Hopf algebras associated to quantum groups are those of prime interest.

We say that a linear functional $\int$ on a left-covariant differential calculus
$\omd$ over a Hopf algebra $\sa$ is {\it left-invariant} if ${(\id_\sa\otimes \int)\dom(\om)}
=(\int \om)1$, for all $\om \in \Om$, where 1 is the unit of $\sa$.

Clearly, the restriction of $\int$ to $\sa$ is a left-invariant linear functional on $\sa$; however, it may be equal to zero on $\sa$ (this is frequently the case).

\proclaim{Theorem 4.3} Let $\int$ be a linear functional on a left-covariant differential calculus $\omd$ over a Hopf algebra $\sa$. Suppose also that $\sa$ admits a Haar integral~$h$. Then the following are equivalent conditions:

(1) $\int a\om = h(a)\int\om$, for all $a\in \sa$ and for all $\om\in \ominv$;

(2) $\int$ is left-invariant. \endproclaim

\demo{Proof} Assume first that $\int$ is left-invariant and suppose that $a\in \sa$ and $\om\in \ominv$. Since $h(1)=1$, we have $\int a\om =h((\int a\om )1)=
{h((\id_\sa \otimes \int)\dom(a\om))}={(h\otimes \int)(\D(a)(1\otimes \om))}
={\int ((h\otimes \id_\sa)\D(a))\om}={\int h(a)\om}=h(a)\!\int\!\om$.
Hence, Condition~(2) implies Condition~(1).

Now suppose that Condition~(1) holds, and let $a$ and $\om$ be as before. We may write $\D(a)={\sum_{i=1}^M b_i\otimes c_i}$, for some elements $b_i$ and $c_i$ in $\sa$. Then
${(\idts)(\dom(a\om))}={(\idts)(\D(a)(1\otimes \om))}=
{(\idts)(\sum_{i=1}^M b_i\otimes c_i\om)}=
{\sum_{i=1}^M (\int c_i\om)b_i}={\sum_{i=1}^M h(c_i)(\int \om)b_i}\allowmathbreak=
{(\id_\sa\otimes h)(\D(a))\int\om}={h(a)(\int\om)\,1}={(\int a\om)\,1}$. Since $\Om$ is the linear span of the elements $a\om$, it follows that $\int$ is left-invariant. Hence, Condition~(1) implies Condition~(2).\qed\enddemo

It is a well-known and useful result that if $h$ is a left-invariant linear functional on a Hopf algebra~$\sa$ and $\k$ is the co-inverse on~$\sa$, then
$$\k((\idth)(\D(a)(1\otimes b)))=(\idth)((1\otimes a)\D(b)),$$
for all elements $a,b\in \sa$.
We show now that a corresponding such result holds for left-invariant linear functionals on a differential calculus.

\proclaim{Theorem 4.4} Let $\omd$ be a left-covariant differential calculus over a Hopf algebra~$\sa$ and let $\int$ be a left-invariant linear functional on~$\Om$. Then,
$$\k((\idts)(\dom(\om)(1\otimes \om')))=(\idts)((1\otimes \om)\dom(\om')),$$
for all $\om,\om'\in \Om$, where $\k$ is the co-inverse of $\sa$. \endproclaim

\demo{Proof} Choose a linear basis $(\om_i)_{i\in I}$ for $\ominv$ and choose linear functionals $f_{ij}$ on $\sa$ such that $\om_ia=\sum_{j\in I} (f_{ij}*a)\om_j$, for all $a\in \sa$; we can do this by Theorem~4.1. Let $E_{ij}$ be the operator on $\sa$ defined by setting $E_{ij}(a)=f_{ij}*a$. It suffices to show the equation in the theorem in the case that $\om=a\om_i$ and $\om'=b\eta$, for arbitrary elements $a,b\in \sa$ and $\eta\in \ominv$ and arbitrary index~$i$.
We have
$$\align
\k &((\id_\sa\otimes \int)(\dom(a\om_i)(1\otimes b\eta)))
=\k((\idts)(\D(a)(1\otimes \om_i)(1\otimes b)(1\otimes \eta)))\\
&=\sum_{j\in I} \k((\idts)(\D(a)(1\otimes E_{ij}(b))(1\otimes \om_j\eta)))=
\sum_{j\in I} \k((\idth_j)(\D(a)(1\otimes E_{ij}(b))),
\endalign$$
where $h_j$ is the linear functional on $\sa$ defined by setting
$h(a)=\int a\om_j\eta$. Since $\int$ is left-invariant, and
$\om_i\eta$ is also left-invariant, so is $h_j$. Hence,
$$\k((\idth_j)(\D(a)(1\otimes E_{ij}(b)))=(\idth_j)((1\otimes a)\D
E_{ij}(b)).$$ However, ${\D E_{ij}(b)}={\D((\id_\sa\otimes
f_{ij})\D(b))}= (\id_\sa\otimes (\id_\sa\otimes f_{ij})\D)\D(b)
=(\id \otimes E_{ij})\D(b)$. \newline Consequently, we have
$$\align
\sum_{j\in I}& \k((\idth_j)(\D(a)(1\otimes E_{ij}(b)))= \sum_{j\in
I} (\idth_j)((1\otimes a)\D E_{ij}(b))\\ &= \sum_{j\in I}
(\idts)((1\otimes a)(\id_\sa\otimes E_{ij})\D(b)(1\otimes \om_j)
(1\otimes \eta))\\ &= (\idts)((1\otimes a)(1\otimes
\om_i)\D(b)(1\otimes \eta))= (\idts)((1\otimes
a\om_i)\dom(b\eta)).\endalign$$ This proves the theorem. \qed
\enddemo

\proclaim{Theorem 4.5} Let $\omd$ be a left-covariant differential calculus over a Hopf algebra~$\sa$ admitting a Haar integral~$h$. Then the linear map, ${P\colon \Om\to \Om}$, defined by setting $P={(h\otimes \id_\Om)\dom}$, is idempotent with image equal to $\ominv$; also, $P(\om_1\om\om_2)=\om_1P(\om)\om_2$, for all $\om \in \Om$ and $\om_1,\om_2\in \ominv$. Moreover, $Pd=dP$.
If $\int$ is a left-invariant linear functional on~$\Om$, then $\int P(\om)=\int\om$, for all $\om\in \Om$. \endproclaim

\demo{Proof} If $a$ is an element of $\sa$ and $\om$ is an invariant form of $\Om$, then $$P(a\om)={(h\otimes \id_\Om)(\D(a)\dom(\om))}=
{(h\otimes \id_\Om)(\D(a)(1\otimes \om))}={(h\otimes \id_\sa)(\D(a))\om}=h(a)\om.$$
It follows from this calculation, that $P(\om)=\om$ and, using the fact that $\Om$ is the linear span of the elements $a\om$, that $P(\Om)=\ominv$. Hence, $P^2=P$.

Now suppose that $\om$ is an arbitrary form of $\Om$ and that $\om_1,\om_2\in \ominv$.
Then $P(\om_1\om\om_2)={(\htid)((1\otimes \om_1)\dom(\om)(1\otimes \om_2))}=
{\om_1(\htid)(\dom(\om))\om_2}={\om_1P(\om)\om_2}$.

We also have $Pd(\om)=
{(\htid)\dom d(\om)}={(\htid)(\id_\sa\!\otimes d)\dom(\om)}=
{d(\htid)\dom(\om)}\allowmathbreak=dP(\om)$. Hence, $Pd=dP$.

Suppose now $\int$ is a left-invariant linear functional on $\Om$. Then $\int \!P(\om)=
{\int \!(h\!\otimes \id_\sa)\!\dom(\om)}\allowmathbreak={(h\otimes \int)\dom(\om)}=
{h((\id_\sa \otimes \int)\dom(\om))}={h((\int \om)1)}=\int\om$. \qed\enddemo

If $\om'$ and $\om$ are invariant elements of $\Om$, then
${\int \om' a \om}={h(a)\int\om'\om}$, since ${\int \om' a \om}=
{\int P(\om' a \om)}={\int \om'P(a) \om}={h(a)\int\om'\om}$.

\proclaim{Corollary 4.6} The linear space of $N$-dimensional, left-invariant linear functionals on $\Om$ is linearly isomorphic to the linear dual of $\ominv_N$. Hence, $\Om$ admits a unique
non-zero, $N$-dimensional, left-invariant linear functional, up to a non-zero scalar factor, if, and only if, $\dim(\ominv_N)=1$. \endproclaim

\demo{Proof} It follows directly from the theorem that he restriction map, ${\int\mapsto \int_{\ominv_N}}$, is the linear isomorphism of the preceding statement. Surjectivity of this map is the only non-obvious point. This is seen by observing that if $\t$ is a linear functional on $\ominv_N$, then we can define the corresponding linear functional on $\Om$ by setting $\int \om=0$, if $\om$ is a $k$-form for which $k<N$, and by setting $\int \om =\t P(\om)$, if $\om\in \Om_N$. Then if $a\in \sa$ and $\om\in \ominv_N$, and if $\D(a)={\sum_{i=1}^M b_i\otimes c_i}$, for some elements $b_i$ and $c_i$ belonging to $\sa$, we have ${(\id_\sa\otimes \int)(\dom(a\om))}=
{(\id_\sa\otimes \int)(\D(a)(1\otimes \om))}={\sum_{i=1}^M \t (P(c_i)\om)b_i}=
{\sum_{i=1}^M h(c_i)\t(\om)b_i}=
{(\idth)(\D(a))\t(\om)}={h(a)\t(\om)1}={\t P(a\om)}={(\int a\om)1}$. Hence, by Theorem~4.1, $\int$ is left-invariant. \qed \enddemo

A Haar integral $h$ on a Hopf algebra $\sa$ is necessarily {\it left faithful} in the sense that, whenever $a$ is an element of $\sa$ for which $h(ba)=0$, for all $b\in \sa$, we must have $a=0$.

\proclaim{Theorem 4.7} Let $\int$ be a non-zero, left-invariant linear functional on a
left-covariant differential calculus~$\omd$ over a Hopf algebra~$\sa$ admitting a Haar integral~$h$.
Then $\int$ is weakly faithful. \endproclaim

\demo{Proof} Suppose that $a\in\sa$ and that $\int \om a =0$, for all $\om\in \Om$. Since $\int\ne 0$, we may choose $\om$ such that $\int \om\ne 0$. Then, for all $b\in \sa$, we have ${0=\int\om ba}={\int P(\om b a)}={\int P(\om)h(ba)}=
{(\int \om)h(ba)}$. It follows, from faithfulness of~$h$, that $a=0$. Hence, $\int$ is weakly faithful. \qed\enddemo

\proclaim{Theorem 4.8} Let $\omd$ be an $N$-dimensional left-covariant differential calculus over the Hopf algebra~$\sa$ admitting a Haar integral~$h$. If $\omd$ admits a left faithful, left-invariant, $N$-dimensional linear functional~$\int$, then $\dim(\ominv_N)=1$. \endproclaim

\demo{Proof} Let $\om$ be an invariant $N$-form of $\Om$ for which
$\int \om =0$. If $a\in \sa$, then $\int a\om=h(a)\int\om =0$. It follows, by faithfulness of $\int$, that $\om=0$. Therefore, the linear map, ${\int\colon \ominv_N\to \bC}$, is injective.
Since $\int$ is non-zero and left invariant, this restriction map cannot be the zero map. Hence, it is a linear isomorphism of $\ominv_N$ onto $\bC$. Therefore, $\dim(\ominv_N)=1$, as required. \qed\enddemo

\proclaim{Corollary 4.9} The functional $\int$ is closed if, and only if, $d(\ominv_{N-1})=0$.
If $\int$ is closed, it is necessarily a twisted graded trace. \endproclaim

\demo{Proof} First observe that if $P=(\htid)\dom$, and $a\in \sa$ and $\om\in \ominv$, then
$\int (da)\om= \int P((da)\om)=\int P(da)\om=\int (dP(a))\om=0$, since $P(a)\in \bC1$ and $d1=0$.
Hence, $\int d(a\om)=\int ad\om+\int (da)\om = \int ad\om$. Using the identification
$\Om_{N-1}=\sa\ominv_{N-1}$, it follows from this observation that if $d(\ominv_{N-1})=0$, then $\int d=0$; that is, $\int$ is closed.
Suppose now conversely that $\int$ is closed and let $\om\in \ominv_{N-1}$. Then $0=\int d(a\om)= \int ad\om$, for all $a\in \sa$. By faithfulness of $\int$, $d(\om)=0$. Hence, $d(\ominv_{N-1})=0$, as required.

Now suppose that $\int$ is closed and we shall show it is a twisted graded trace. Choose any
non-zero element $\ta$ in $\ominv_N$ for which $\int \ta=1$; then $\ominv_N=\bC\ta$. Since $\sa\ta=\ta\sa$, by Theorem~4.1, there is a unique automorphism $\r_1$ of $\sa$ such that $\ta a=\r_1(a)\ta$, for all $a\in \sa$. Also, the Haar integral $h$ admits another automorphism $\r_2$ of $\sa$ such that $h(ba)=h(\r_2(a)b)$, for all $a,b\in \sa$. Set $\s_0=\r_2\r_1$. Then
$\int b\ta a=\int b\r_1(a)\ta = h(b\r_1(a)) = h(\r_2\r_1(a)b)=\int \s_0(a)b\ta$. It follows from Theorem~1.2 that $\int$ is a twisted graded trace. \qed\enddemo

We say that an $N$-dimensional differential calculus $\omd$ over a unital algebra $\sa$ is {\it non-degenerate} if, whenever $\om$ is a $k$-form in $\Om$ for which $\om'\om=0$, for all $\om'\in \Om_{N-k}$, we necessarily have $\om=0$. It is clear that if $\Om$ admits a left faithful,
$N$-dimensional linear functional, then $\Om$ is non-degenerate.

\proclaim{Theorem 4.10} Let $\omd$ be a non-degenerate, $N$-dimensional, left-covariant differential calculus over a Hopf algebra $\sa$ admitting a Haar integral~$h$. Then $\Om$ admits a left faithful, left-invariant,
$N$-dimensional linear functional if, and only if, $\dim(\ominv_N)=1$. \endproclaim

\demo{Proof} The forward implication follows from Theorem~4.6.
Suppose conversely  $\dim(\ominv_N)=1$. Then, by Theorem~4.6,
$\Om$ admits a non-zero, $N$-dimensional, left-invariant linear
functional $\int$ (unique up to multiplication by a scalar
factor). To prove the theorem, we have only to show now that
$\int$ is left faithful. Thus, we must show that if $\om\in \Om$ and
$\int \om'\om=0$, for all $\om'\in \Om$, then $\om=0$. We may
clearly suppose, without loss of generality, that $\om\in \Om_k$,
for some index $k\le N$. Then if $\om'\in \Om_{N-k}$, we have
$\om'\om=a\ta$, for some element $a\in \sa$. Hence, if $b\in \sa$,
$\int b\om'\om=0$, by assumption. Hence, $h(ba)=0$, for all $b\in
\sa$. By faithfulness of~$h$, $a=0$. Therefore, $\om'\om=0$. We
now use non-degeneracy of $\Om$ to deduce that $\om=0$, as
required. \qed\enddemo

Woronowicz has constructed a certain non-degenerate, left-covariant, three-dimensional calculus $\omd$ over the Hopf algebra $\sa$ underlying the compact quantum group $SU_q(2)$, where $q$ is a real parameter for which $0<\vert q\vert\le 1$. For this calculus, $\ominv_1$ has a linear basis $\om_0,\om_1,\om_2$ for which $\sa\om_i=\om_i\sa$, for ${i=0,1,2}$. Hence, for each index~$i$, there exists an automorphism $\r_i$ of $\sa$ such that $\om_ia=\r_i(a)\om_i$, for all $a\in \sa$.

Since $SU_q(2)$ is a compact quantum group, it admits a Haar integral~$h$. Also, there is an automorphism $\r$ of $\sa$ such that $h(ba)=h(\r(a)b)$, for all $a,b\in \sa$. We define
a 1-dimensional, left-invariant linear functional $\int$ on $\Om$ by setting
${\int a_0\om_0+a_1\om_1+a_2\om_2}=h(a_1)+h(a_2)$. This functional is closed. To see this, observe first that there exist linear functionals $\chi_0,\chi_1,\chi_2$ on $\sa$ such that $da=\sum_{i=0}^2 (\chi_i*a) \om_i$, for all $a\in \sa$. Since ${d1=0}$, we have $\chi_i(1)=0$, for all~$i$. Using this, and right-invariance of $h$, we get $\int da=
{h(\chi_1*a)+h(\chi_2*a)}={h(a)\chi_1(1)+h(a)\chi_2(1)=0}$.

We claim now that $\int$ is not a twisted graded trace. Otherwise, let $\s$ denote its twist automorphism. Then ${\int \om a}={\int \s(a)\om}$, for all $\om\in \Om_1$. Therefore, for $a,b\in \sa$ and $i=1,2$, we have $h(\r\r_i(a)b)=
h(b\r_i(a))={\int b\r_i(a)\om_i}={\int b\om_i a}={\int \s(a)b\om_i}=h(\s(a)b)$. Faithfulness of $h$ now implies that $\r\r_i(a)=\s(a)$, for all $a\in \sa$. Hence, $\r_1=\r_2$. But if $\a,\g$ are the canonical generators of $SU_q(2)$ as in~\cite{\WORA}, then $\r_1(\a)=q^{-2}\a$ and
$\r_2(\a)=q^{-1}\a$, by Table~1 of~\cite{\WORA}. Hence, $\r_1\ne \r_2$. This contradiction shows that, as claimed, $\int$ is not a twisted graded trace.

We now truncate Woronowicz's calculus to get a 1-dimensional differential calculus
${(\Om',d')}$ over $\sa$. Then ${(\Om',d')}$ is a non-degenerate, left-covariant, one-dimensional calculus over~$\sa$, and $\om_0,\om_1,\om_2$ is a linear basis for the space of invariant 1-forms.

The restriction $\int'$ of $\int$ to $\Om'$ is a closed, left-invariant, 1-dimensional linear functional on $\Om'$. As we saw is the case for $\int$, the functional $\int'$ is also not a twisted graded trace. This shows that the faithfulness hypothesis in Theorem~4.10 is necessary.

\proclaim{Lemma 4.11} Let $\int$ be a left-invariant twisted graded trace on the universal unital differential calculus $\tomd$ over a Hopf algebra $\sa$ admitting a Haar integral~$h$. Let $I$ be the left kernel of~$\int$ and let $J={I\cap \tilde{\Om}^{\text{inv}}}$. Then the linear map from ${\sa\otimes J}$ to $I$ that sends ${a\otimes \om}$ onto  $a\om$ is an isomorphism of left $\sa$-modules.
Hence, $I$ is invariant under $\D_{\tilde{\Om}}$ in the sense that $\D_{\tilde{\Om}}(I) \subseteq {\sa \otimes I}$.
\endproclaim

\demo{Proof} Let $\om\in I$; using the identification of ${\sa\otimes \tilde{\Om}^{\text{inv}}}$ with $\tilde{\Om}$, we write, as we may, $\om={\sum_{i=1}^M a_i\om_i}$, where ${a_1,\dots,a_M}$ are linearly independent elements of $\sa$, and ${\om_1,\dots\om_M}$ are left-invariant elements of $\tilde{\Om}$. Set ${X = \{(h(b a_1),\ldots,h(b a_M)) \mid b \in \sa\}}$. We claim that $X = \bC^M$. Suppose otherwise (and we shall obtain a contradiction). Then there exists a non-zero linear functional $\t$ on $\bC^M$ such that $\t(x) = 0$, for all $x \in X$. Clearly, $\t$ is determined by scalars $\mu_1,\ldots,\mu_M$, in the sense that ${\t(\lambda_1,\ldots,\lambda_M)} = {\sum_{i=1}^M \lambda_i\mu_i}$, for all ${\lambda_1,\ldots,\lambda_M \in \bC}$. Moreover, since $\t\ne 0$, the scalars $\mu_i$ are not all equal to zero. Now let $b\in \sa$. Then
${h(b(\sum_{i=1}^M \mu_ia_i))}={\sum_{i=1}^M \mu_i h(b a_i)} =
{\t(h(b a_1),\ldots,h(b a_M)) =0}$.
Hence, ${\sum_{i=1}^M \mu_ia_i = 0}$, by faithfulness of $h$. This contradicts the linear independence of the elements $a_1,\ldots,a_M$.
Consequently, to avoid contradiction, we must have $X = \bC^M$. It follows that there exist elements ${b_1,\ldots,b_M \in \sa}$ such that ${h(b_j a_i) = \delta_{ji}}$,
for ${i,j=1,\ldots,M}$. Hence, for any invariant element $\eta$ in $\tilde{\Om}$, we have, since $\om\in I$, ${0=\sum_{i=1}^M\int \eta b_ja_i\om_i}={\sum_{i=1}^M h(b_ja_i)\int \eta\om_i}={\int \eta\om_j}$. Therefore, for any element $a\in \sa$, ${\int a\eta \om_j}={h(a)\int \eta\om_j=0}$. Consequently, the form $\om_j$ belongs to $I$ and therefore, since it is left-invariant, it belongs to $J$. The lemma now follows. \qed\enddemo

\proclaim{Theorem 4.12} Let $\bint$ be an $N$-dimensional, left-invariant, closed twisted graded trace on the universal unital differential calculus $\tomd$ over a Hopf algebra $\sa$ admitting a Haar integral~$h$.
The $N$-dimensional differential calculus $\omd$ associated to $(\tom,d,\bint)$ is left-covariant and the canonical twisted graded trace $\int$ on $\omd$ is left-invariant. \endproclaim

\demo{Proof} Let $\phi$ be the canonical algebra isomorphism from ${A\otimes \Om}$ onto the quotient algebra ${(\sa\otimes \tom)/(\sa\otimes I)}$ obtained by mapping ${a\otimes (\om+I)}$ onto
${a\otimes \om+\sa\otimes I}$, for all $a\in \sa$ and $\om\in \tom$. Then the map,
${\dom \colon \Om\to \sa\otimes \Om}$, defined by setting $\dom(\om+I)=
{\phi^{-1}(\D_{\tilde{\Om}}\,\om+\sa\otimes I)}$ for all $\om\in \tilde{\Om}$, is a co-action making $\omd$
left-covariant. This follows from the readily verified facts that $\dom$ is an algebra homomorphism extending the co-multiplication on $\sa$ and that ${(\id_\sa\otimes d)\dom}=\dom d$.

To see that $\int$ is left-invariant, let $\om\in \tom$ and
suppose that $\D_{\tilde{\Om}}(\om)= {\sum_{i=1}^Ma_i\otimes \om_i}$, for
some elements $a_i$ in $\sa$ and forms $\om_i$ in $\bom$. Then
$$\align (\id_\sa\otimes \int)\dom(&\om+I)=(\idts)(\sum_{i=1}^M
a_i\otimes (\om_i+I))= \sum_{i=1}^M (\int \om_i+I)a_i\\
&=\sum_{i=1}^M (\bint \om_i)a_i=(\id_\sa\otimes \bint)(\D_{\tilde{\Om}}\, \om)
=(\bint \om)1=(\int\om+I)1.\endalign$$ Thus, $\int$ is
left-invariant, as required. \qed\enddemo

Let $\sa$ be a Hopf algebra.
Let $\bar \sa$ be the quotient linear space $\sa/\bC 1$ and, for $a\in \sa$,
write $\bar a$ for ${a+\bC}$ in $\bar\sa$. We define a map $\D_N$ from ${\sa\otimes \bar\sa^{\otimes N}}$ to ${\sa\otimes (\sa\otimes \bar\sa^{\otimes N})}$ by setting ${\D_N(a_0\otimes \bar a_1\otimes \cdots \otimes \bar a_N)}=
{(m\otimes \id)(FG(a_0\otimes \bar a_1\otimes \cdots \otimes \bar a_N))}$. Here, $m$ is the unique linear map from ${\sa^{\otimes (N+1)}}$ to $\sa$ that sends the
elementary tensor ${a_0\otimes \cdots\otimes a_N}$ onto the product ${a_0\cdots a_N}$. The map $F$ is the
unique linear map from ${(\sa\otimes \sa)\otimes (\sa\otimes \bar\sa)^{\otimes N}}$ to ${\sa^{\otimes (N+2)}\otimes \bar\sa^{\otimes N}}$ that maps the elementary tensor ${(b_0\otimes
c_0)\otimes (b_1\otimes \bar c_1)\otimes \cdots \otimes (b_N\otimes \bar c_N)}$ onto ${(b_0\otimes \cdots \otimes
b_N)\otimes (c_0\otimes \bar c_1\cdots \otimes \bar c_N)}$. Finally,
$G$ is the $(N\!+\!1)$-fold tensor product ${\D\otimes \bar\D \otimes \cdots \otimes \bar\D}$, where $\bar\D$ is the algebra homomorphism from $\bar\sa$ to $\sa \otimes \bar\sa$ defined as follows. Let $\pi$ be the quotient map from $\sa$ to $\bar\sa$, then $\bar\D$ is defined by setting $\bar\D(\pi(a)) = (\id_\sa \otimes \pi)\D(a)$ for all $a \in \sa$.

Straightforward calculations show that $\D_N$ is a left co-action of $\sa$ on
${\sa\otimes \bar\sa^{\otimes N}}$.

If ${\f\colon \sa^{N+1}\to \bC}$ is a multilinear function that vanishes on any element ${(a_0,a_1,\dots,a_N)}$, whenever any of the components ${a_1,\dots,a_N}$ belongs to $\bold C$1, we let $\hat \f$ be the corresponding linear map on
${\sa\otimes \bar\sa^{\otimes N}}$
(so that ${\hat\f(a_0\otimes  \bar a_1\cdots \otimes \bar a_N)}= {\f(a_0,\dots,a_N)}$).
We say that $\f$ is {\it left-invariant} if
${(\id_\sa \otimes \hat\f)\D_N(c)}={\hat\f(c)1}$, for all $c\in
{\sa\otimes \bar \sa^{\otimes N}}$, where 1 is the unit of $\sa$.

Suppose now that $\f$ is the twisted cyclic cocycle associated an
$N$-dimensional, closed twisted graded trace $\int$ on $\Om$, for some
left-covariant differential calculus $\omd$ over $\sa$. A straightforward calculation shows that
$$(\id\otimes \int)(\D(a_0)\dom d(a_1)\cdots \dom d(a_N))=(\id\otimes \hat\f)(\D_N(a_0\otimes \bar a_1\otimes \cdots  \otimes \bar a_N)),$$
for all elements ${a_0,a_1,\dots,a_N\in \sa}$. From this it follows easily that $\int$ is left-invariant if, and only if, $\hat\f$ is left-invariant.

We summarise our observations here in the following result.

\proclaim{Theorem 4.13} Suppose that $\omd$ is a left-covariant differential calculus over a Hopf algebra $\sa$ and that $\int$ is an $N$-dimensional closed, twisted graded trace on $\Om$. Let $\f$ be the corresponding twisted cyclic $N$-cocycle. Then $\int$ is left-invariant if, and only if, $\f$ is left-invariant. \endproclaim

\subheading{5. A construction of a 3-dimensional differential calculus}

In this section we show how our construction of a differential calculus from a closed twisted graded trace on the universal unital differential calculus can be used to show the existence of a 3-dimensional calculus first constructed by very different means by Woronowicz.

First, recall that the universal unital differential calculus $\tom$ over a Hopf~algebra $\sa$ is left covariant. Let $\kappa$ be the co-inverse on~$\sa$, and denote by $m$ the linear map from ${\sa\otimes \tom}$ to $\tom$ that sends the elementary tensor ${a\otimes \om}$ onto the product $a\om$. Define the linear map $w$ from $\sa$ to $\tom^{\text{inv}}_1$ by setting $w(a)={m(\kappa\otimes d)\D(a)}$. If the unit 1 of $\sa$ and the family $(e_i)_{i\in I}$ form a linear basis for $\sa$, then, for each positive integer~$k$, the products of the form
$w(e_{i_1})\cdots w(e_{i_k})$, where ${i_1,\dots,i_k\in I}$, form a linear basis of $\tom^{\text{inv}}_k$, that we shall call the basis {\it associated to} the family~$(e_i)$~\cite{\KS}.

If $\sa$ is a Hopf $*$-algebra, then $\tom$ is a $*$-differential calculus over $\sa$, where $w(a)^*= -w(\kappa(a)^*)$, for all $a\in \sa$. Here, as usual, $\kappa$ is the co-inverse of $\sa$.

Suppose now that $q$ is a non-zero real parameter for which ${\vert q \vert \le 1}$. We denote by $\sa_q$ the Hopf algebra associated to the compact quantum group~$SU_q(2)$~\cite{\WORA}. Recall that $\sa_q$ is the universal unital $*$-algebra generated by a pair of elements~$\a$ and $\g$ satisfying the relations
$$\align
\a^*\a+\g^*\g=1 &\qquad \a\a^*+q^2\g\g^*=1\\
\g^*\g=\g\g^* \qquad \a\g&=q\g\a\qquad \a\g^*=q\g^*\a.
\endalign$$
The co-multiplication $\D$ on $\sa_q$ is the unique unital $*$-homomorphism for which
$\D(\a)={\a\otimes \a}-{q\g^*\otimes \g}$ and $\D(\g)={\g\otimes \a}+
{\a^*\otimes \g}$.

Let $\bold E=\bZ\times \bN\times \bN$. For $\e=(k,l,m)\in \bold E$, denote by $a_\e$ the product $\a^k \gamma^l(\gamma^*)^m$. Here we use the usual convention in this context that for $k<0$, $\a^k=(\a^*)^{-k}$. It is well known that these elements $a_\e$ form a linear basis for $\sa_q$, that we shall call the {\it standard} basis. Writing $w_\e$ for $w(a_\e)$, it follows that the non-zero products $w_{\e_1}w_{\e_2}w_{\e_3}$ form a basis for $\tom^{\text{inv}}_3$, that we shall call the {\it standard} basis of $\tom^{\text{inv}}_3$
(of course, the element $w_{\e_1}w_{\e_2}w_{\e_3}$ is equal to zero only if $\e_i=(0,0,0)$, for some index~$i$).

Again suppose that $\e=(k,l,m)$. We set $c(\e)=0$ if $l$ or $m$ are positive and we set $c(\e)=c(k) =(1-q^{-2k})(1-q^{-2})^{-1}$, if $l=m=0$. If $\om$ is a standard basis element, $\om=w_{\e_1}w_{\e_2}w_{\e_3}$, we set $c(\om)=
{c(\e_1)+c(\e_2)+c(\e_3)}$.

We shall say that $\e$ is {\it reduced}, or that $a_\e$ is reduced, if $(k,l,m) \not= (0,0,0)$ and if $(k,l)=(0,1)$, $(0,0)$ or $(1,0)$; in this case we set $t(\e)=-1$, 0, or $1$, respectively, and we call $t(\e)$ the {\it type} of $\e$.

We shall say that a standard basis element $\om=w_{\e_1}w_{\e_2}w_{\e_3}$ is {\it reduced}, if all the factors have reduced indices and their types are distinct. We set
$t(\om)=(t(\e_1),t(\e_2),t(\e_3))$, and call this triple the {\it type} of $\om$.

Using Theorem~4.3, we define a 3-dimensional left-invariant linear functional
$\int$ on the universal unital differential calculus $\tom$ over $\sa_q$ by setting
$\int$ equal to zero on all of the non-reduced standard basis elements, and by defining $\int$ on a reduced standard basis element $\om=w_{\e_1}w_{\e_2}w_{\e_3}$ as follows:

\vbox{(1) if $t(\om)=(-1,0,1)$, $\int \om= c(\om)$; \ \ \ \ \ \ \ \ \ (4) if $t(\om)=(0,1,-1)$, $\int \om= q^6c(\om)$;

(2) if $t(\om)=(-1,1,0)$, $\int \om= -q^4c(\om)$; \ \ \ \ (5) if $t(\om)=(1,-1,0)$, $\int \om=  q^6c(\om)$;

(3) if $t(\om)=(0,-1,1)$, $\int \om= -q^4c(\om)$; \ \ \ \ (6) if $t(\om)=(1,0,-1)$, $\int \om= -q^{10}c(\om)$.}

Henceforth, we denote by $e$ the co-unit of $\sa_q$. Since $e(\a)=1$ and $e(\g)=0$, we have $e(a_\e)=0$, unless $\e=(k,0,0)$, for some $k\in \bZ$, in which case $e(a_\e)=1$.
If $a\in \sa_q$, we write $\bar a$ for $a-e(a)1$ in the sequel.

\proclaim{Theorem 5.1} The functional $\int$ is self-adjoint. \endproclaim

\demo{Proof} If $\om$ belongs to $\tom$, we have to show that $\int \om^*=(\int\om)^-$.
Using linearity and left-invariance of $\int$, we may clearly reduce to the case where $\om=w_{\e_1}w_{\e_2}w_{\e_3}$, for some elements ${\e_1,\e_2,\e_3\in \bold E}$.
It follows from the formulas $w(a)^*=-w(\kappa(a)^*)$ and ${\k(a(k,l,m))^*}=
{(-1)^{l+m}q^{l-m}a(k,m,l)}$ that we have
${w_{(k,l,m)}^*}={(-1)^{l+m+1}q^{l-m}w_{(k,m,l)}}$.
Hence, $\om^*=-w_{\e_3}^* w_{\e_2}^* w_{\e_3}^*$ is a scalar times a standard basis element that is reduced if, and only if, $\om$ is reduced. Therefore, $\int \om^*=0=(\int \om)^-$, if $\om$ is non-reduced.

We now suppose that $\om$ is reduced. In this case each index $\e_i$ is reduced and we have ${\sum_{i=1}^3 \e_i(2)=\sum_{i=1}^3 \e(3)=1}$, where $\e_i=(\e_i(1),\e_i(2),\e_i(3))$.
Write ${\e'_i=(\e_i(1),\e_i(3),\e_i(2))}$. Then $\om^*=
{(-1)^{\sum_i \e_i(2)+\e_i(3)}q^{\sum_i \e_i(3)-\e_i(2)} w_{\e'_3}w_{\e'_2} w_{\e'_1}}
={w_{\e'_3}w_{\e'_2} w_{\e'_1}}$. We have only to show now, since $\int \om$ is
real-valued, that $\int \om =\int \om'$, where $\om'={w_{\e'_3}w_{\e'_2} w_{\e'_1}}$.
Now if the type of $\om$ is ${t(\om)=(t_1,t_2,t_3)}$, then the type of $\om'$ is clearly
${t(\om')=(-t_3,-t_2,-t_1)}$. An examination of the definition of $\int$ given in Equations~(1)--(6) above shows that $\int$ takes on the same values on reduced standard basis elements if one is of type ${(t_1,t_2,t_3)}$ and the other of type
${(-t_3,-t_2,-t_1)}$. Hence, $\int \om =\int \om'$ and the proof is completed. \qed\enddemo

\proclaim{Lemma 5.2} Let $a,b,c \in \sa_q$ and suppose that $c=\g$ or $c=\g^*$. Suppose also that $\om_1,\om_2\in \tom^{\text{inv}}_1$. Then
$\int a  w(\bar{b} c)  \om_1
\om_2 = \int a  \om_1 w(\bar{b} c)   \om_2 = \int a \om_1
\om_2\,  w(\bar{b} c) = 0$. \endproclaim

\demo{Proof} We shall prove only that the first integral vanishes and only in the case that $c=\g$. The other cases have similar proofs. We may clearly suppose, by left-invariance of $\int$, that $a=1$. By linearity of $\int$, we may also suppose that $\om_1=w_{\e_1}$ and $\om_2=w_{\e_2}$ and that $b=a_\e$, for some elements $\e,\e_1,\e_2\in \bold E$.

If $e(b)=0$, then $b\g=a_{\e+(0,1,0)}$ is not reduced. Hence, the element $w(\bar bc)\om_1\om_2$ is not reduced. It follows that $\int$ vanishes on this element.

On the other hand, if $e(b)=1$, then $\bar b={a_\e-e(a_\e)}={\a^k-1}$, for some $k\in \bZ$. Hence, ${\int w(\bar bc)\om_1\om_2}={\int w_{(k,1,0)}w_{\e_1}w_{\e_2}}-
{\int w_{(0,1,0)}w_{\e_1}w_{\e_2}}$. This is equal to zero, since it follows easily from the definition of $\int$ that ${\int w_{(k,1,0)}w_{\e_1}w_{\e_2}}$ is independent of~$k$. \qed\enddemo

Let $\s_0$ be the twist automorphism associated to the Haar measure~$h$ on~$\sa_q$; that is, $\s_0$ is the unique automorphism on $\sa_q$ for which $h(a'a)=h(\s_0(a)a')$, for all $a,a'\in \sa_q$. Let $\s_1$ be the unique automorphism on~$\sa_q$ for which
$\sigma_1(\al) = {q^{-4} \al}$, ${\sigma_1(\gamma) = q^{-4}  \gamma}$, ${\sigma_1(\al^*) = q^4 \al^*}$ and ${\sigma_1(\gamma^*) = q^4  \gamma^*}$. (This automorphism exists as a consequence of the universal property enjoyed by~$\sa_q$). Finally, set
$\sigma = \sigma_0 \sigma_1$; of course, $\s$ is again an automorphism.

\proclaim{Lemma 5.3} For all $a\in \sa_q$ and $\om \in \tom_3$, we have
${\int \omega a} = {\int \sigma(a) \omega}$. \endproclaim

\demo{Proof} To see this, we may use the fact that $\s$ is an automorphism to reduce to the case where $a$ is equal to one of ${\a,\g,\a^*}$ or $\g^*$. We may also use simple linearity of~$\int$ to reduce to the case that $\om=b\om'$, where $b\in \sa_q$ and $\om'=w_{\e_1}w_{\e_2}w_{\e_3}$, for some elements ${\e_1,\e_2,\e_3\in \bold E}$, and ${\om'\ne 0}$.
We write $\e_i={(k_i,l_i,m_i)}$. We shall consider in this proof only the case that $a=\a$; the other three cases have similar proofs.

We begin by observing that it is easy to check that there exists a positive integer~$M$ and elements ${c_{10},\dots,c_{1M}}$, $\dots$, ${c_{40},\dots,c_{4M}}$ in $\sa_q$ such that
$${(\id_\sa\ot \id_\sa\ot \D)(\id_\sa \ot \D)\D(\al)} =
{\sum_{i=0}^M c_{1i} \ot c_{2i} \ot c_{3i} \ot c_{4i}}$$
and such that $c_{j0}=\a$, for ${j=1,2,3,4}$ and either $\gamma$ or $\gamma^*$ occurs in the triple ${c_{2i},c_{3i},c_{4i}}$, for each index ${i=1,\dots,M}$. It follows from~\cite{\KS, 14.3.2~Eqn.~(51)} that
${w(a_{\e_1})w(a_{\e_2})w(a_{\e_3})\a}={\sum_{i=0}^M c_{1i}w(\bar a_2 c_{2i})w(\bar a_3 c_{3i})w(\bar a_4 c_{4i})}$,
for all $a_2,a_3,a_4\in \sa_q$.
This, and Lemma~5.2, gives
$$\align   \int b &w(a_{\e_1})w(a_{\e_2})w(a_{\e_3})  \al  \\
&= \int b  \al w(\overline{a_{\e_1}}\al)w(\overline{a_{\e_2}}\al)w(\overline{a_{\e_3}} \al) +\sum_{i=1}^n \int b  c_{1i}w(\overline{a_{\e_1}}c_{2i})w(\overline{a_{\e_2}}c_{3i})
w(\overline{a_{\e_3}} c_{4i})\\
&=  \int b\al w(\overline{a_{\e_1}}\al)w(\overline{a_{\e_2}}\al)w(\overline{a_{\e_3}} \al). \endalign$$

We now divide the proof up into two cases: where $\om'$ is reduced and where it is non-reduced. Considering first the case that $\om'$ is non-reduced, we have $\e_i$ is
non-reduced, for some index~$i$, so that $l_i+m_i \ge 2$. Hence,
${w(\overline{a_{\e_i}} \al)}={w(a_{\e_i} \al)}
 = {\l w(a_{(k_i+1,l_i,m_i)})} + {\mu  w(a_{(k_i+1,l_i+1,m_i+1)})}$,
for some scalars $\l,\mu \in \bC$. Clearly, ${a_{(k_i+1,l_i,m_i)}}$ and
${a_{(k_i+1,l_i+1,m_i+1)}}$ are non-reduced. It follows from this that
${\int \om \a}=
{\int \!b\al w(\overline{a_{\e_1}}\al)w(\overline{a_{\e_2}}\al)w(\overline{a_{\e_3}} \al)=0}$.
Also, by left-invariance of $\int$, we have ${\int \s(\a)\om}=
{h(\s(\a)b)\int\om'=0}$. Thus, ${\int \om\a}={\int \s(\a)\om}$, if $\om'$ is non-reduced.

We turn now to the second case, where $\om'$ is reduced. First observe that for any integer~$k$, we have
$$\align
&   w(\overline{\al^k \gamma}  \al)  =  w(\al^k \gamma  \al) =
q^{-1}  w(\al^{k+1} \gamma) + \l  w(\al^{k+1} \gamma^2 \gamma^*) \\
&  w(\overline{\al^k \gamma^*}  \al)  =  w(\al^k \gamma^*  \al) =
q^{-1}  w(\al^{k+1} \gamma^*) + \mu  w(\al^{k+1} \gamma (\gamma^*)^2) \tag{$*$}\\
&  w(\overline{\al^k }  \al) = w((\al^k -1)
\al) = w(\al^{k+1} )  -  w(\al) + \nu  w(\al^{k+1}\gamma \gamma^*),
\endalign$$
for some scalars $\l,\mu,\nu \in \bC$ (precise knowledge of these scalars is
not needed).

We now have six subcases to consider, depending on the type $t(\om')$ of $\om'$.
We shall only give the details for the case that $t(\om')={(0,-1,1)}$---the other five cases are proved by similar means.

Using Equations~($*$), we get
$$\align
\int b  w(\al^{k_1})  w(\al^{k_2} \gamma^*)  w(\al^{k_3}
 \gamma)  \al & =  q^{-2}  \int b \al
(w(\al^{k_1+1})-w(\al))  w(\al^{k_2+1} \gamma^*)
w(\al^{k_3+1} \gamma)  \\  & =  q^{-2}  h(b \al)  (-q^4) (c(k_1+1)-1).
\endalign$$
One can easily check that $c(k_1+1) -1 = q^{-2} c(k_1)$. Hence,
$$\align
\int \om\a &=  \int  b  w(\al^{k_1})  w(\al^{k_2} \gamma^*)  w(\al^{k_3}\gamma)\al
= q^{-4}h(\sigma_0(\al) b)   (-q^4) \,c(k_1) \\
&= q^{-4}   \int \sigma_0(\al) b  w(\al^{k_1})  w(\al^{k_2}\gamma^*)  w(\al^{k_3}\gamma) = \int \s(\a)\om.
\endalign$$
This completes the proof. \qed\enddemo

We turn now to the problem of showing that $\int$ is closed. To this end, we introduce some convenient notation. If ${m\colon \tom\ot\tom \to \tom}$ is the multiplication map, we write $\tilde w(c)$ for ${m(w\ot w)(c)}$, for all
${c\in \sa_q\ot \sa_q}$. If $a\in \sa_q$, then
$d(w(a)) = \tilde{w}(\D(a))$~\cite{\KS, 14.3.2~Eqn.~(52)}.

It is easy to verify that, for all $\om_1,\om_2,\om_3\in \tom^{\text{inv}}_1$, we have
${\int \om_1\om_2\om_3=0}$, if any one of the factors $\om_i$ is of the form $w_\e=w(a_\e)$, where $\e$ is a non-reduced element of $\bold E$.

Let $J$  be the left ideal of ${\sa_q \ot \sa_q}$ generated by the elements
${\gamma^2 \ot 1}$, ${\gamma^{*2} \ot 1}$, ${\gamma \gamma^* \ot 1}$, ${1 \ot
\gamma^2}$, ${1 \ot \gamma^{*2}}$ and ${1 \ot \gamma \gamma^*}$. Using the commutation relations in the axioms for $\sa_q$, one easily sees that $J$ is also a right ideal of ${\sa_q\ot \sa_q}$. On $\sa_q \ot \sa_q$ we denote by $\cong$ the equivalence relation determined by $J$, so that $c \cong c'$ if, and only if, $c'-c \in J$. Using linearity of $\int$ and the fact that $\int$ vanishes on all non-reduced standard basis elements of~$\tom_3$, one can easily verify that, for all ${c,c'\in \sa_q\ot\sa_q}$ such that $c\cong c'$, and all 1-forms~$\om$, ${\int \om \tilde{w}(c)} = {\int \om \tilde{w}(c')}$ and ${\int \tilde{w}(c) \om} = {\int \tilde{w}(c') \om}$.

We shall make frequent, often tacit, use of these observations in the sequel of this section.

\proclaim{Theorem 5.4} The linear functional $\int$ is closed. \endproclaim

\demo{Proof} By left-invariance, it is clear that we need only show that
${\int d\om=0}$, if $\om\in \tom^{\text{inv}}$. Using the fact that $\int$ is non-zero only on the 3-forms, we may further suppose that $\om\in \tom^{\text{inv}}_2$. Using linearity of $\int$ and $d$, we may even suppose that $\om=\om_{\e_1}\om_{\e_2}$, for some elements $\e_1,\e_2\in \bold E$.
Given these considerations, the result will now follow from the next four lemmas.
\qed\enddemo

\proclaim{Lemma 5.5} Let $a,b \in \sa_q$ and suppose that $b$ is equal to one of
$\gamma^2$, $\gamma^{*2}$ or $\gamma \gamma^*$. Let $\om \in \tom^{\text{inv}}_1$. Then
$\int d(w(ab)\om) = \int d(\om w(ab)) = 0$. \endproclaim

\demo{Proof} We may, and shall, suppose that $\om=w_{\e}$, for some $\e\in \bold E$.
We shall show only that $\int d(w(ab)\om) = 0$; the proof that ${\int d(\om w(ab)) = 0}$ is similar.

First, observe that ${\int w(ab)\om'=0}$, for any $\om'\in \tom^{\text{inv}}_2$. (To see this, one can reduce to the case that $a=a_{\e_1}$ and $\om'=w_{\e_2}w_{\e_3}$.) Hence,
${\int d(w(a b) \om)} = {\int d(w(ab))\om} - {\int w(ab)d\om}={\int d(w(ab))\om}$.

We consider first the case where $b =  \gamma^2$. Then
$$\align
\de(ab)  & =  \de(a)(\gamma \ot \al + \al^* \ot \gamma)^2  \\
& =\de(a)(\gamma^2 \ot \al^2 + \gamma \al^* \ot \al \gamma + \al^* \gamma \ot
\gamma \al + \al^{*2} \ot \gamma^2) \\
& \cong  \de(a)( q \al^* \gamma
\ot \al \gamma + q^{-1} \al^* \gamma \ot \al \gamma). \endalign$$
Hence, $\int d(w(ab))  \om =  \int \tilde{w}(\de(a)(q \al^*\ot \al+q^{-1}
\al^* \ot \al)(\g\ot\g))  \om = 0$,
where the last equality follows from the easily verified fact that
${\int \tilde w(c(\g\ot\g))\om=0}$, for all ${c\in \sa_q\ot \sa_q}$ (to see this, one can take ${c=a_{\e_1}\ot a_{\e_2}})$.

The second case, where $b=\g^{*2}$, is dealt with similarly.

We turn now to the third case, where $b = \gamma \gamma^*$. Using the first two cases already proved, we may suppose now, without loss of generality, that
$a = \al^p$, for some $p \in \bZ$. Then
$$\align
\de(ab) & =  \de(a)(\gamma \ot \al + \al^* \ot \gamma)  (\gamma^* \ot
\al^* + \al \ot \gamma^*) \\
&  =  \de(a)(\gamma \gamma^* \ot \al \al^* +
\gamma \al  \ot \al \gamma^* + \al^* \gamma^* \ot \gamma \al^* + \al^* \al
\ot \gamma \gamma^*) \\
& \cong  \de(\al^p) ( q^{-1} \al \gamma \ot \al
\gamma^* + q  \al^* \gamma^* \ot \al^* \gamma).\endalign$$
Since $\de(\al) = \al \ot \al - q \gamma^* \ot \gamma$, the above
chain of equivalences implies that
$$\align \de(ab) & \cong (\al^p \ot \al^p) ( q^{-1} \al \gamma \ot \al \gamma^* + q
\al^* \gamma^* \ot \al^* \gamma)\\
& \cong  ( q^{-1} \al^{p+1} \gamma \ot \al^{p+1} \gamma^* + q
\al^{p-1} \gamma^* \ot \al^{p-1} \gamma).
\endalign$$
Therefore,
$$\int d(w(a b))  \om  =  q^{-1}  \int w(\al^{p+1}
\gamma)w(\al^{p+1} \gamma^*)  \om + q  \int w(\al^{p-1}
\gamma^*)w(\al^{p-1} \gamma)  \om. \tag $*$ $$
Write $\e = (k,l,m)$.
If $l+m \ge 1$, both terms on the right hand side of Equation~($*$) are equal to zero.
If $l=m=0$, we have ${\int d(w(ab))  \om} = {q^{-1} (q^6 c(k)) + q (-q^4 c(k)) = 0}$.
Thus, ${\int d(w(a b))  \om =0}$, whatever the value of $\e$. \qed\enddemo

\proclaim{Lemma 5.6} Let $a,b,c \in \sa_q$ and suppose that $c$ is equal to either $\g$ or $\g^*$.  Then ${\int d(w(a c)w(b c)) = 0}$. \endproclaim

\demo{Proof} We shall prove the result only in the case that $c= \gamma$; the case that $c=\g^*$ has a similar proof. First, we remark that ${\int \om_1\om_2\om_3=0}$, if the elements $\om_1,\om_2,\om_3\in \tom^{\text{inv}}_1$, and two of them are of the form $w(f\g)$ and $w(f'\g)$, for some elements $f,f'\in \sa_q$. To see this, we may use linearity to reduce to the case where the other factor is of the form $w_{\e''}$, and $f=a_\e$ and $f'=a_{\e'}$. Then $\om_1\om_2\om_3$ will be either the zero element, or a standard basis element that is not reduced, so that ${\int \om_1\om_2\om_3=0}$. We use this now to show ${\int d(w(a c)w(b c)) = 0}$.

We may write $\de(a)= {\sum_{i=1}^M a^1_i \ot a^2_i}$ and
$\de(b) = {\sum_{i=1}^M b^1_i \ot b^2_i}$, for some elements ${a^1_i,a^2_i,b^1_i,b^2_i\in \sa_q}$. Then, since $\de(\gamma) = \gamma \ot \al + \al^* \ot \gamma$,
we have
$$\align
\int d(w(a \gamma) w(b \gamma))  & =  \int \tilde{w}(\de(a \gamma))
 w(b \gamma) - \int w(a \gamma)  \tilde{w}(\D(b \gamma)) \\
& = \sum_{i=1}^M \int w(a^1_i \gamma) w(a^2_i \al) w(b \gamma) +
    \sum_{i=1}^M \int w(a^1_i \al^*) w(a^2_i \gamma) w(b \gamma) \\
& \ \ \ - \sum_{i=1}^M \int w(a \gamma)w(b^1_i \gamma)  w(b^2_i \al)
  -\sum_{i=1}^M \int w(a \gamma)w(b^1_i \al^*)  w(b^2_i \gamma).
\endalign$$
It follows from the remarks in the preceding paragraph that all the terms in the four sums vanish, and therefore ${\int d(w(a \gamma) w(b \gamma)) = 0}$, as required. \qed\enddemo

\proclaim{Lemma 5.7} Let $k_1$ and $k_2$ be integers and suppose that $c$ is equal to  $\g$ or $\g^*$. Then ${\int d(w(\al^{k_1}c) w(\al^{k_2} c^*)) = 0}$. \endproclaim

\demo{Proof} The two cases have similar proofs; we give the proof for $c=\g$ only. We have
$$\align \de(\al^{k_1} \gamma) & =  \de(\al^{k_1})(\gamma \ot \al + \al^* \ot
\gamma)  \cong (\al^{k_1} \ot \al^{k_1})( \gamma \ot \al + \al^* \ot \gamma) \\
&\cong  \al^{k_1} \gamma \ot \al^{{k_1}+1} + \al^{{k_1}-1} \ot \al^{k_1} \gamma; \tag{5.1} \endalign$$
and similarly,
$$\de(\al^{k_2} \gamma^*) \cong  \al^{k_2} \gamma^* \ot \al^{{k_2}-1} +
\al^{{k_2}+1} \ot \al^{k_2} \gamma^*.$$
Therefore,
$$\align \int d(w(\al^{k_1} \gamma)w(\al^{k_2} \gamma^*))
& =  \int d(w(\al^{k_1} \gamma))  w(\al^{k_2}\gamma^*)  - \int w(\al^{k_1} \gamma) d(w(\al^{k_2} \gamma^*)) \\
& =  \int w(\al^{k_1} \gamma) w(\al^{{k_1}+1}) w(\al^{k_2} \gamma^*) + \int w(\al^{{k_1}-1})w(\al^{k_1} \gamma) w(\al^{k_2} \gamma^*) \\
&\ - \int w(\al^{k_1} \gamma)w(\al^{k_2} \gamma^*) w(\al^{{k_2}-1}) - \int w(\al^{k_1} \gamma) w(\al^{{k_2}+1})
w(\al^{k_2} \gamma^*) \\
& =  - q^{10} c({k_1}+1) + q^6 c({k_1}-1) - q^6 c({k_2}-1) +q^{10} c({k_2}+1). \endalign$$
An elementary calculation shows that ${c(k-1)- q^4  c(k+1)} = {-q^2(1+q^2)}$,
for all ${k \in \bZ}$, and it follows that ${\int d(w(\al^{k_1} \gamma)  w(\al^{k_2}
\gamma^*)) = 0}$. \qed\enddemo

\proclaim{Lemma 5.8} If $k_1$ and $k_2$ are integers and $c$ is equal to 1, $\g$ or $\g^*$, then ${\int w(\a^{k_1}c)w(\a^{k_2})}={\int w(\a^{k_1})w(\a^{k_2}c)=0}$. \endproclaim

\demo{Proof} We show the result only in the case $c=1$ and $c=\g$; the case $c=\g^*$ is proved similarly to that for $c=\g$. We shall assume that $k_1$ and $k_2$ are positive; the cases where they are not can be dealt with by similar arguments to that we present here.

Let $m$ be a positive integer.
Since $c(1) = 1$, $c(m+1) = 1 + q^{-2}c(m)$, and
$(\gamma^* \ot \gamma)(\al \ot \al) = q^{-2} \, (\al \ot \al)(\gamma^* \ot \gamma)$, a simple induction argument
based on the formula $\de(\al) = \al \ot \al - q \gamma^* \ot \gamma$ implies that ${\de(\al^m)} =
{(\al \ot \al- q \gamma^* \ot \gamma)^m} \cong {\al^m \ot \al^m-qc(m)(\al^{m-1} \ot \al^{m-1})(\gamma^* \ot \gamma)}$.
Consequently,
$$\align
\int d(w(\al^{k_1})w(\al^{k_2}))
& =  \int w(\al^{k_1}) w(\al^{k_1})w(\al^{k_2})-qc({k_1})
\int w(\al^{{k_1}-1} \gamma^*) w(\al^{{k_1}-1}\gamma) w(\al^{k_2}) \\
& \ \ \ - \int w(\al^{k_1})  w(\al^{k_2}) w(\al^{k_2})
+ q  c({k_2}) \int w(\al^{k_1})  w(\al^{{k_2}-1}\gamma^*)w(\al^{{k_2}-1} \gamma) \\
& =  0 + q^5\, c({k_1}) c({k_2}) -0 - q^5  c({k_2}) c({k_1}) =0.
\endalign$$

Now we consider the case where $c=\g$.
By Equation~(5.1) in the proof of Lemma~5.7, we know that
${\de(\al^{k_1} \gamma)} \cong {\al^{k_1} \gamma \ot \al^{{k_1}+1}}
+ {\al^{{k_1}-1} \ot\al^{k_1}\gamma}$. By the considerations of the preceding paragraph, we know that
${\de(\al^{k_2})} \cong {\al^{k_2} \ot \al^{k_2}}
+ {\lambda\,\al^{{k_2}-1} \gamma^* \ot\al^{{k_2}-1} \gamma}$, for some $\lambda \in \bC$. Hence,
$$\align \int d(w(\al^{k_1} \gamma) w(\al^{k_2}))
& = \int w(\al^{k_1} \gamma)w(\al^{{k_1}+1}) w(\al^{k_2}) + \int w(\al^{{k_1}-1}) w(\al^{k_1} \gamma)w(\al^{k_2}) \\
&  \ - \int w(\al^{k_1} \gamma) w(\al^{k_2})w(\al^{k_2})
- \lambda \int w(\al^{k_1} \gamma) w(\al^{{k_2}-1}\gamma^*) w(\al^{{k_2}-1}\gamma). \endalign$$
All of the terms in the sum are easily seen to be zero, since they involve evaluating $\int$ on non-reduced standard basis elements. Therefore,
${\int d(w(\al^{k_1} \gamma)w(\al^{k_2})) =0}$. \qed\enddemo

Now that we have established that $\int$  is closed, we use Lemma~5.3 and Theorem~1.1 to deduce that $\int$ is a twisted graded trace. Moreover, the twist automorphism $\tilde{\s}$ of $\int$ extends the automorphism~$\s$ of $\sa_q$. We use these facts and the fact that $\int$ is self-adjoint, to apply the construction of Section~1 to the triple $(\tilde{\Om},d,\int)$ to deduce the existence of a left-covariant,
3-dimensional $*$-differential calculus $\Om$ over $\sa_q$. We shall denote the canonical twisted graded trace on $\Om$ by the same symbol $\int$ and refer to the domains of these functionals to distinguish them in cases of ambiguity.

Let $\pi$ denote the quotient map from $\tilde{\Om}$ onto $\Om$.
It is easy to verify from the definition of $\int$ on $\tilde{\Om}$ that,

(1) For all $k\in \bZ$, $\pi(w(\al^k)) = c(k)\pi(w(\al))$, $\pi(w(\al^k \gamma)) =\pi(w(\gamma))$ and $\pi(w(\al^k \gamma^*)) =  \pi(w(\gamma^*))$;

(2) For all $k,l,m\in \bZ$ for which $l,m\ge 0$ and $l+m\ge 2$, we have ${\pi(w_{(k,l,m)}) = 0}$.

Set $\om_0 = - q  \pi(w(\gamma^*))$, $\om_1= \pi(w(\al))$ and $\om_2 =
- q^{-1} \pi(w(\gamma))$. It follows from Conditions~(1) and~(2) that $\om_0$, $\om_1$ and $\om_2$ linearly span $\ominv$.
It is immediate from the definition of $\int$ on $\tilde{\Om}$ that
$$\matrix
 \int \om_0  \om_1  \om_2 & = & 1 & \ \  &
\int \om_0  \om_2  \om_1 & = & - q^4
\\ \int \om_1  \om_0  \om_2 & = & -q^4   &  \ \   &
\int \om_1  \om_2  \om_0 & = &  q^6 & \qquad\qquad (5.2)\\
\int \om_2  \om_0  \om_1 &
= & q^6 &  \ \  & \int \om_2  \om_1  \om_0 & = & - q^{10} \endmatrix$$
and that $\int \om_i \om_j \om_k = 0$ for every $i,j,k \in \{0,1,2\}$
where any two of the indices $i,j,k$ are the same.

Since the trace $\int$ on $\Om$ is left faithful, it follows easily that $\om_0$, $\om_1$ and $\om_2$ are linearly independent and therefore that they form a
linear basis for $\ominv$.

Let $a$ and $b_1,\dots,b_M$ and $c_1,\dots,c_M$ be elements in $\sa_q$ such
that $\de(a) = {\sum_{i=1}^M b_i \ot c_i}$.  Then, by Equations~(51) and (52) of \cite{\KS, 14.3.2} , and the equation $w(a)^*= -w(\kappa(a)^*)$, which holds for all $a\in \sa_q$, we have

(1) $\pi(w(a))b =  {\sum_{i=1}^M  b_i  \pi(w(\overline{a}c_i))}$, for all
$b \in \sa_q$;

(2) $d a = \sum_{i=1}^M b_i\, \pi(w(c_i))$;

(3) $d\pi(w(a)) = \sum_{i=1}^M \pi(w(b_i))\,\pi(w(c_i))$;

(4) $\om_0^*=q\om_2\qquad \om_1^*=-\om_1\qquad \om_2^*=q^{-1}\om_0$.

\smallskip
Applying these formulas in our particular case, it is easy to check that
the differential calculus $(\Om,d)$ that we have constructed here
satisfies the formulas in Tables~1, 2 and~6 of~\cite{\WORA}.

Using faithfulness $\int$ on $\Om$, combined with the formulas
in Equations~(5.2), one can readily verify that our differential calculus also satisfies the formulas of Table 5 of~\cite{\WORA} and that the three
elements $\om_0\om_1$, $\om_0\om_2$ and $\om_1\om_2$ form a linear basis for
$\ominv_2$.

With this information at hand, it is now straightforward to conclude that
our $*$-differential calculus $(\Om,d)$ is isomorphic to the 3-dimensional calculus constructed by Woronowicz in~\cite{\WORA} by an entirely different method.

We believe that our method for constructing calculi is one that is perhaps more natural
than other methods, since the basis of our approach is essentially to devise a  ``presentation'' of the calculus in terms of generators and relations. We hope that the method will be used in the future to construct other interesting calculi.

\enddocument

\Refs

\ref \no  1 \by  E. Abe \book Hopf Algebras \publ Cambridge University Press
\publaddr Cambridge \yr 1980 \endref

\ref \no 2 \by A. Connes \book Noncommutative Geometry \publ Academic Press
\publaddr San Diego--New York \yr 1994 \endref

\ref \no 3 \by C. Kassel \book Quantum Groups \publ Springer \publaddr Berlin--Heidelberg \yr 1995 \endref

\ref \no 4 \by A.U. Klymik and K. Schm\"udgen \book Quantum Groups and their Representations \publ Springer \publaddr Heidelberg--New York \yr 1998 \endref

\ref \no 5 \by J. Kustermans and L. Tuset \paper A survey of C*-algebraic quantum groups~I \jour Irish Math. Soc. Bull. \vol 43 \yr 1999 \pages 8--63 \endref

\ref \no 6 \by G.J. Murphy and L. Tuset \paper Compact quantum groups, preprint, National University of Ireland, Cork \yr 1999 \endref


\ref \no 7 \by S.L. Woronowicz \paper Twisted $SU(2)$ group---an example of
a non-commutative differential calculus \jour Publ. RIMS Kyoto Univ. \vol 23 \yr 1987 \pages 117--181 \endref

\ref \no 8 \by S.L. Woronowicz \paper Compact matrix pseudogroups
\jour Comm. Math. Phys. \vol 111 \yr 1987 \pages 613--665 \endref

\ref \no 9 \by S.L. Woronowicz, \paper Differential calculus on compact matrix pseudogroups (quantum groups) \jour Commun. Math. Phys. \vol 122 \yr 1989 \pages 125--170 \endref

\endRefs
\bye